\documentclass[aos]{imsart}
\usepackage{amsmath}
\usepackage[colorlinks=true]{hyperref}
\usepackage{amsfonts}
\usepackage{epsfig}
\usepackage{amssymb}
\usepackage{graphicx}

\def\cT{{\cal T}}
\def\cN{{\cal N}}

\def\mysix{8}

\def\Riskopt{\Risk_{\hbox{\tiny\rm opt}}}
\def\RiskS{{\hbox{\rm RiskS}}}
\def\RiskoptS{\hbox{\rm RiskS}_{\hbox{\tiny\rm opt}}}
\def\tRiskoptS{\hbox{\rm \scriptsize RiskS}_{\hbox{\tiny\rm opt}}}

\def\cl{\mathop{\hbox{\rm cl}}}
\def\inter{\mathop{\hbox{\rm int}}}

\def\cond{\mathop{\hbox{\rm\small Cond}}}
\def\I{{\cal I}}

\def\brA{{{A}}}
\def\brB{{{B}}}
\def\brX{{\Y}}
\def\brn{{\bar{n}}}

\def\bT{{\mathbf{T}}}
\def\bK{{\mathbf{K}}}
\newcommand{\half}{ \mbox{\small$\frac{1}{2}$}}
\newcommand{\four}{\mbox{\small$\frac{1}{4}$}}

\def\Risk{{\hbox{\rm Risk}}}
\def\Ker{\mathop{\hbox{\rm Ker}}}

\def\Z{{\cal Z}}

\def\Diag{{\hbox{\rm Diag}}}
\def\bR{{\mathbf{R}}}

\def\Opt{{\hbox{\rm Opt}}}

\def\A{{\cal A}}
\def\L{{\cal L}}
\def\U{{\cal U}}

\def\T{{\cal T}}

\def\Q{{\cal Q}}

\def\X{{\cal X}}
\def\Y{{\cal Y}}
\def\N{{\cal N}}

\def\S{{\cal S}}

\def\bE{{\mathbf{E}}}
\def\bS{{\mathbf{S}}}
\def\Prob{\hbox{\rm Prob}}

\newtheorem{lemma}{Lemma}[section]
\newtheorem{corollary}{Corollary}[section]
\newtheorem{proposition}{Proposition}[section]

\newtheorem{theorem}{Theorem}[section]
\def\Det{{\hbox{\rm Det}}}
\def\qed{$\Box$}

\def\e{{\rm e}}
\def\Tr{{\hbox{\rm Tr}}}

\def\tRiskopt{\hbox{\rm\scriptsize Risk}_{\hbox{\rm\tiny opt}}}
\newcommand{\aic}[2]{{\color{blue}~#2}}

\newcommand{\be}{\begin{eqnarray}}
\newcommand{\ee}[1]{\label{#1}\end{eqnarray}}
\newcommand{\nn}{\nonumber \\}
\newcommand{\ese}{\end{eqnarray*}}
\newcommand{\bse}{\begin{eqnarray*}}
\newcommand{\rf}[1]{~(\ref{#1})}
\newcommand{\wh}[1]{{\widehat{#1}}}
\newcommand{\hide}[1]{{}}
\begin{document}
\begin{frontmatter}
\title{Near-Optimality of Linear Recovery
 in Gaussian Observation Scheme under $\|\cdot\|_2^2$-Loss}
\runtitle{Near-optimality of linear recovery}
\begin{aug}
\author{\fnms{Anatoli} \snm{Juditsky} \thanksref{t1}\ead[label=e1]{anatoli.juditsky@univ-grenoble-alpes.fr}}
\and
\author{\fnms{Arkadi} \snm{Nemirovski} \thanksref{t2}\ead[label=e2]{nemirovs@isye.gatech.edu}}
\thankstext{t1}{This author was supported by the  LabEx PERSYVAL-Lab (ANR-11-LABX-0025)}
\thankstext{t2}{Research of this author was supported by NSF grants   CCF-1523768 and CMMI-1262063.}
\runauthor{A. Juditsky and A. Nemirovski}
\affiliation{Universit\'e Grenoble Alpes and Georgia Institute
 of Technology}
\address{LJK, Universit\'e Grenoble-Alpes\\ 700 Avenue Centrale\\
38401 Domaine Universitaire de Saint-Martin-d'H\`eres\\ France
\printead{e1}}
\address{Georgia Institute
 of Technology\\ Atlanta, Georgia
30332\\ USA\\
\printead{e2}}
\end{aug}
\begin{abstract}
We consider the problem of recovering linear image $Bx$ of a signal $x$ known to belong to a given convex compact set $\X$ from indirect observation
$\omega=Ax+\sigma\xi$ of $x$ corrupted by Gaussian noise $\xi$. It is shown that under some assumptions on $\X$ (satisfied, e.g., when $\X$ is the intersection of $K$ concentric ellipsoids/elliptic cylinders),
an easy-to-compute linear estimate is near-optimal in terms of its worst-case, over $x\in \X$, expected $\|\cdot\|_2^2$-loss. The main novelty here is that the {result} imposes no restrictions
on $A$ and $B$.  To the best of our knowledge, preceding results on optimality of linear estimates dealt either with one-dimensional $Bx$ (estimation of linear forms) or with the  ``diagonal case'' where $A$, $B$ are diagonal and ${\X}$ is given by a ``separable'' constraint like ${\X}=\{x:\sum_ia_i^2x_i^2\leq 1\}$ or ${\X}=\{x:\max_i|a_ix_i|\leq1\}$.
\end{abstract}


\end{frontmatter}
\section{Introduction}
In this paper we address one of the most basic problems of High-Dimensional Statistics, specifically, as follows: given indirect noisy observation
 $$
 \omega=Ax+\sigma\xi\eqno{[A:m\times n, \;\xi\sim \N(0,I_m)]}
 $$
 of unknown ``signal'' $x$ known to belong to a given convex compact subset ${\X}$ of $\bR^n$, we want to recover the image $w=Bx\in\bR^\nu$ of $x$ under a given linear mapping. We focus on the case where the quality of a candidate recovery $\omega\mapsto{\wh{w}}(\omega)$ is quantified by its worst-case, over $x\in {\X}$, expected $\|\cdot\|_2^2$-error, that is, by the risk
 $$
 \Risk[{\wh{w}}|{\X}]=\sup_{x\in {\X}}\left[{\bE_\xi}\left\{ \|{\wh{w}}(Ax+\sigma\xi)-Bx\|_2^2\right\}\right]^{1/2}.
 $$
 The simplest and the most studied type of recovery is affine one: ${\wh{w}}(\omega)=H^T\omega+h$; assuming ${\X}$ symmetric w.r.t. the origin, we lose nothing when passing from affine estimates to linear ones -- those of the form ${\wh{w}}_H(\omega)=H^T\omega$.
{Starting from the pioneering works \cite{kuks1,kuks2}, linear estimates received much attention in the statistical literature (see, e.g., \cite{Rao1973,Rao1976,Cris,Dru,Pilz,Arnold} and references therein).} An advantage of linear estimates, from the computational point of view, is that under favorable circumstances (e.g., when ${\X}$ is an ellipsoid), minimizing risk over linear estimates is an efficiently solvable problem.
{ On the other hand, linear estimates are also of major importance to statistical theory.
For instance, a huge body of results on {\em rate-optimality} of linear estimates on various signal classes (which arise from some classes of regular functions) form the backbone of classical nonparametric statistics (see, e.g., \cite{IbrHas1981,wasserman2006all,Tsybakov}).
Furthermore, for several important signal classes linear estimates occur to be near-optimal on the class of {\em all} possible estimates. This is, for instance, the case for signal recovery from direct observations (the case of $B=A=I$) in the situation where the set $\X$ of signals is an ellipsoid or a box. The case of ellipsoidal $\X$ was studied first by M.S. Pinsker, see \cite{Pinsker1980}, who showed that in the problem of recovery of the signal $x\in {\X}$ from direct observation $\omega=x+\sigma\xi$, ${\X}$ being a ``Sobolev ellipsoid'' of the form $\{x\in\bR^n:\,\sum_jj^{2\alpha}x_j^2\leq L^2\}$,  the ratio of the risk of a properly selected linear estimate to the minimax risk $\Riskopt[{\X}]:=\inf_{{\wh{w}}(\cdot)}\Risk[{\wh{w}}|{\X}]$  (the infimum is taken over all estimates, not necessarily linear) tends to 1, as $\sigma\to +0$, and this happens {\sl uniformly in $n$}, $\alpha$ and $L$ being fixed. Similar ``asymptotic optimality'' results are also known {for} ellipsoids related to classes of analytic functions \cite{golubev1996asymptotically} and for ``diagonal'' case, where ${\X}$ is the above ellipsoid/box and $A$, $B$ are diagonal matrices  \cite{efromovich1996sharp} (see also \cite{efromovich2008nonparametric} for modern presentation of that approach).  {The results on non-asymptotic near-optimality of linear estimates (up to a factor 1.11...) are also available for the case where $A=B=I$ and $\X$ is an ellipsoid ($\X=\{x\in\bR^n:\,\sum_ja_j^{2}x_j^2\leq 1\}$ for given $a_j$) or a box ($\X=\{x\in\bR^n:\,\max_j|a_jx_j|\leq 1\}$) (see, e.g., \cite{donoho1990minimax}) (the corresponding argument can be easily extended to the case of diagonal $A$ and $B$).} Note that the situation is quite different for the problem of estimation of a {\em linear form} $w=b^Tx$ (i.e., the case of one-dimensional $Bx$).
An exceptional from several points of {view} ``general'' (that is, not imposing severe restrictions on how the geometries of ${\X}$, $A$ and $B$ are linked to each other) result on optimality of linear estimates in this case is due to D. Donoho who proved \cite{Don95} that when recovering a linear form, the best, over all linear estimates, risk is within the factor 1.11... of the minimax risk.
}
 \par
The goal of this paper is to establish a rather general result on near-optimality of properly built linear estimates as compared to all possible estimates. Note that a result of this type is bounded to impose some restrictions on ${\X}$, since there are cases (e.g., the one of a high-dimensional $\|\cdot\|_1$-ball ${\X}$) where linear estimates are {\sl by far} suboptimal. Our restrictions on the family of sets (we call them {\em ellitopes}) ${\X}$ reduce to the existence of a special type representation of ${\X}$ and are satisfied, e.g., when ${\X}$ is the intersection of $K<\infty$ ellipsoids/elliptic cylinders:
  \begin{equation}\label{exampleX}
  {\X}=\{x\in\bR^n: x^TS_kx\leq 1,1\leq k\leq K\}\quad\quad[S_k\succeq0,\sum_kS_k\succ0].
 \end{equation}
In particular, ${\X}$ can be a symmetric w.r.t. the origin compact polytope given by $2K$ linear inequalities $-1\leq s_k^Tx\leq 1$, $1\leq k\leq K$. Another instructive example is a set of the form ${\X}=\{x:\|Sx\|_p\leq L\}$, where $p\geq2$ and $S$ is a matrix with trivial kernel. It should be stressed that while imposing some restrictions on ${\X}$, {\sl we  require nothing from $A$ and $B$}. Our main result (Theorem \ref{prop1}) states, in particular, that in the case of ${\X}$ given by (\ref{exampleX}) and arbitrary $A$, $B$, the risk of properly selected linear estimate ${\wh{w}}_{H_*}$, with both $H_*$ and the risk being efficiently computable, satisfies the bound
 \[
 \Risk[{\wh{w}}_{H_*}|{\X}]\leq O(1)\sqrt{\ln\left({O(1)\|B\|^2K^2\kappa^{-1}\over\tRiskopt^2[{\X}]}\right)}\Riskopt[{\X}],\eqno{(*)}
 \]
 where $\|B\|$ is the spectral norm of $B$, $\kappa$ is the minimal eigenvalue of $\sum_kS_k$, $\Riskopt[{\X}]$ is the minimax risk, and $O(1)$ stands for an absolute constant. {It should be mentioned that technique used to construct lower bound for optimal risks leads to more precise oracle inequalities when imposing constraints on the structure of the signal class $\X$ and matrices $A,\; B$; in particular, it allows to reproduce classical ``asymptotic'' optimality  results, e.g., in the situation considered in \cite{Pinsker1980,efromovich1996sharp}. On the other hand, we do not know if the bound $(*)$ can be significantly improved in some important ``simple cases'', for instance, in the case where
 $B=I$ and $\X$ is an ellipsoid, without imposing any restrictions on $A$. In this work, however, we prefer to see our approach as} ``operational'' -- {the {\sl provably nearly optimal} estimate itself, its risk and even the lower risk bound involved are all} given by an efficient computation { which is supposed to provide precise near-optimality guaranties for each set of the problem data. From this point of view, the oracle inequality ${(*)}$
 can be viewed as a general indication of a ``goodness'' of linear estimates in a certain context, namely, where the signal set is an intersection of ``not too large'' number of ellipsoids/elliptic cylinders.} This is in sharp contrast with traditional results of non-parametric statistics, where near-optimal estimates and their risks are given in a ``closed analytical form,'' at the price of severe restrictions on the structure of the ``data'' ${\X}$, $A$ and $B$. \footnote{Since this paper has been submitted, the proposed approach has been further developed in \cite{arxiv2017New}. For instance, it is shown that similar near-optimality guaranties for linear estimators can be obtained for more general risks (e.g., $\ell_p$-loss with $1\leq p\leq 2$) and slightly more general sets $\X$, which are solution sets of systems of quadratic {\em matrix} inequalities, for deterministic bounded noises, etc.}
 This being said, it should be stressed  that one of the crucial components of our construction is {completely} classical -- this is the idea, going back to \cite{Pinsker1980} to bound from below the minimax risk via Bayesian risk associated with properly selected Gaussian prior\footnote{\cite{Pinsker1980} addresses the problem of $\|\cdot\|_2$-recovery of a signal $x$ from direct observations ($A=B=I$) in the case where ${\X}$ is a high-dimensional ellipsoid with ``regularly decreasing half-axes,'' like ${\X}=\{x\in\bR^n:\sum_jj^{2\alpha}x_j^2\leq L^2\}$ with $\alpha>0$. In this case
 Pinsker's construction shows that as $\sigma\to +0$, the risk of properly built linear estimate is, uniformly in $n$,  $(1+o(1))$ times the minimax risk. This is much stronger than $(*)$, and it seems quite unlikely that a similarly strong result may hold true in the general case underlying $(*)$.}.
 \par
 The main body of the paper is organized as follows. Section \ref{sect1} contains problem formulation (section \ref{sitgoal}), construction of the linear estimate we deal with (section \ref{sectest}) and the central result on near-optimality of this estimate (section \ref{mainobs}).  Section \ref{sectextensions}  contains some extensions.
 Specifically, we present a version of our main result for the case when the usual worst-case expected $\|\cdot\|_2^2$-risk is replaced with properly defined {\sl relative} risk (section \ref{estrelscale}) and provide a robust, w.r.t. uncertainty in $A,B$, version of the estimate (section \ref{sectrob}). In section \ref{secbypro} we show that the key argument underlying the proof of our main result can be used beyond the scope of statistics, specifically, when quantifying the approximation ratio of the semidefinite relaxation bound on the maximum of a quadratic form over an ellitope.

Proofs are relegated to Appendix.

\section{Situation and main result}
\label{sect1}
\subsection{Situation and goal}\label{sitgoal}
Given $\nu\times n$ matrix $B$, consider the problem of estimating linear image $Bx$ of unknown  signal $x$ known to belong to a given set ${\X}
\subset\bR^n$ via noisy observation
\begin{equation}
\label{eq1obs}
\omega=Ax+\sigma\xi
\end{equation}
where an $m\times n$ matrix $A$ and $\sigma{>}0$ are given, and $\xi\sim\N(0,I_m)$ is the standard Gaussian observation noise.  From now on we assume that ${\X\subset \bR^n}$ is a set given by
\begin{equation}\label{setX}
\X=\big\{x\in\bR^n:\,\exists (y\in\bR^{\brn}, t\in\T):\; x=Py, \,  y^TS_ky\leq t_k,\,1\leq k\leq K\big\},
\end{equation}
{where}
\begin{itemize}
\item $P$ is an $n\times \bar{n}$ matrix,
\item $S_k\succeq0$ are $\bar{n}\times \bar{n}$ matrices with $\sum_kS_k\succ 0$,
\item $\T$ is a nonempty computationally tractable\footnote{for all practical purposes, it suffices to assume that $\T$ is given by an explicit {\sl semidefinite representation}
 $$
 \T=\{t:\exists w: A(t,w)\succeq0\},
 $$
 where $A(t,w)$ is a symmetric and affine in $t,w$ matrix.} convex compact subset of $\bR^K_+$ intersecting the interior of $\bR^K_+$ and such that $\T$ is monotone, meaning that the relations $0\leq\tau\leq t$ and $t\in\T$ imply that $\tau\in\T$.\footnote{The latter relation is ``for free'' -- given a nonempty convex compact set $\T\subset\bR^K_+$, the right hand side of (\ref{setX}) remains intact when passing from $\T$ to its ``monotone hull'' $\{\tau{\in \bR^K_+:\,}\exists t\in\T:\tau\leq t\}$ which already is a monotone convex compact set.} Note that under our assumptions $\inter \T\neq\emptyset$.
 \end{itemize}
 We assume that $BP\neq0$, since otherwise one has $Bx=0$ for all $x\in \X$ and the estimation problem is trivial.
 In the sequel, we refer to a set of the form (\ref{setX}) with data $[P,\{S_k, 1\leq k\leq K\},\T]$ satisfying just formulated assumptions as to an {\sl {ellitope}},
 and to (\ref{setX}) -- as to {\sl {ellitopic} representation} of ${\X}$. Here are instructive examples of {ellitopes} (in all these examples, $P$ is the identity mapping):
\begin{itemize}
\item when $K=1$, $\T=[0,1]$ and $S_1\succ0$, ${\X}$ is the ellipsoid $\{x:x^TS_1x\leq1\}$;
\item when $K\geq1$, $\T=\{t\in\bR^K:\,0\leq t_k\leq 1,\,k\leq K\}$, and ${\X}$ is the intersection
\[\bigcap\limits_{1\leq k\leq K}\{x:\;x^TS_kx\leq 1\}\]
of centered at the origin ellipsoids/elliptic cylinders. In particular, when $U$ is a $K\times n$ matrix of rank $n$ with rows $u_k^T$, $1\leq k\leq K$, and $S_k=u_ku_k^T$, ${\X}$ is symmetric w.r.t. the origin polytope
    $\{x: \|Ux\|_\infty\leq 1\}$;
\item when $U$, $u_k$ and $S_k$ are as in the latter example and $\T=\{t\in\bR^K_+: \,\sum_kt_k^{p/2}\leq1\}$ for some $p\geq 2$, we get
    ${\X}=\{x:\,\|Ux\|_p\leq1\}$.
\end{itemize}
{It should be added that the family of {ellitope}-representable sets is quite rich: this family admits a ``calculus'', so that more {ellitopes} can be constructed by taking intersections, direct products, linear images (direct and inverse) or arithmetic sums of ``basic {ellitopes}'' {given by} the above examples. In fact, the property to be an {ellitope} is preserved by all basic operations with sets preserving convexity and symmetry w.r.t. the origin, see {Section B}.
 \par
 As another instructive, in the context of non-parametric statistics, example of an {ellitope}, consider the situation where our signals $x$ are discretizations of functions of continuous argument running through a compact $d$-dimensional domain $D$, and the functions $f$ we are interested in are those satisfying a Sobolev-type smoothness constraint -- an upper bound on the  $L_p(D)$-norm of $\L f$, where $\L$ is a linear differential operator with constant coefficients. After discretization, this
restriction can be modeled as $\|L x\|_p\leq 1$, with properly selected matrix $L$. As we already know from the above example, when $p\geq2$, the set ${\X}=\{x:\|Lx\|_p\leq1\}$ is an {ellitope}, and as such is captured by our machinery. Note also that by the outlined calculus, imposing on the functions $f$ in question {\em several Sobolev-type smoothness constraints} with parameters $p\geq2$, still results in a set of signals which is an ellitope.}

\paragraph{Estimates and their risks} In the outlined  situation, a candidate estimate is a Borel function ${\wh{w}}(\cdot):\bR^m\to\bR^\nu$; given observation (\ref{eq1obs}), we recover $w=Bx$ as ${\wh{w}}(\omega)$. In the sequel, we quantify the quality of an estimate by its worst-case, over $x\in {\X}$, expected $\|\cdot\|_2^2$ recovery error:
\[
\Risk[{\wh{w}}|{\X}]=\sup_{x\in {\X}}\Big[\bE_{\xi\sim\N(0,I_m)}\left\{\|{\wh{w}}(Ax+\sigma\xi)-Bx\|_2^2\right\}\Big]^{1/2}
\]
and define the optimal, or the {\sl minimax}, risk as
\[
\Riskopt[{\X}]=\inf\limits_{{\wh{w}}(\cdot)}\Risk[{\wh{w}}|{\X}],
\]
where $\inf$ is taken over all Borel candidate estimates.
\par
Our initial observation is that when replacing matrices $A$ and $B$ with $AP$ and $BP$, respectively, we pass from the initial estimation problem of interest -- one where the signal set ${\X}$ is given by \rf{setX}, and we want to recover $Bx$, $x\in {\X}$, via observation (\ref{eq1obs}), to the {\sl transformed problem}, where the signal set is
\[
\brX=\{y\in \bR^{\brn}: \exists t\in \T: \;y^TS_ky\leq t_k,\;1\leq k\leq K\},
\]
and we want to recover $[BP] y$, $y\in \brX$, via observation
\[%
\omega=[AP] y+\sigma\xi.
\]%
It is obvious that the considered families of estimates (the family of all linear and the family of all estimates), same as the risks of the estimates, remain intact under this transformation; in particular,
\[
\Risk[{\wh{w}}|{\X}]=\sup_{y\in\brX}\left[\bE_{\xi}\{\|{\wh{w}}([AP]\,y+\sigma\xi)-[BP]\,y\|_2^2\}\right]^{1/2}.
\]
Therefore, to save notation, from now on and unless mentioned otherwise, we assume that matrix $P$ is identity, so that $\X$ is the ellitope
\begin{equation}\label{Xequals}
\X=\left\{x\in \bR^n:\;\exists t\in \T, \;x^TS_kx\leq t_k,\;1\leq k\leq K\right\}.
\end{equation}
\par\noindent
{\sl Main goal} of what follows is to demonstrate that a  {\sl linear in $\omega$} estimate
\begin{equation}\label{est}
{\wh{w}}_H(\omega)=H^T\omega
\end{equation}
with properly selected efficiently computable matrix $H$ is near-optimal in terms of its risk. We start with building this estimate.
\subsection{Building linear estimate}\label{sectest}
Restricting ourselves to linear estimates (\ref{est}), we may be interested in the estimate with the smallest risk, that is, associated with a $\nu\times m$ matrix $H$ which is an optimal solution to the optimization problem
$$
\min_H\left\{ R(H):=\Risk^2[\wh{w}_H|{\X}]\right\}
$$
We have
\bse
R(H)&=&\max_{x\in \X}\bE_\xi\|H^T\omega-Bx\|^2_2=\bE_\xi\|H^T\xi\|_2^2+\max_{x\in \X}\|H^TAx-Bx\|_2^2\\
&=&\sigma^2\Tr(H^TH)+\max_{x\in \X}x^T(H^TA-B)^T(H^TA-B)x.
\ese
As the maximum over $x$ of convex quadratic functions of $H$, $R(H)$ is itself convex. However, as the maximum over $\X$ of a quadratic in $x$ function, $R(H)$ is typically hard to compute\footnote{For instance, when $\X$ is a unit cube $\{x\in\bR^n:\|x\|_\infty\leq1\}$, computing $R(0)$ in the case of general-type $B$ is equivalent  to maximizing over $X$ a general-type convex quadratic form; it is known that solving the latter problem already within 4\% accuracy  is NP-hard.}. For this reason, we use a linear estimate yielded by minimizing an {\sl efficiently computable convex upper bound} on $R(H)$ which is built as follows.
Let $\phi_\T$ be the support function of $\T$:
$$
\phi_\T(\lambda)=\max_{t\in \T}\lambda^Tt: \bR^K\to\bR.
$$
Observe that whenever $\lambda\in\bR^K_+$ and $H$ are such that
\begin{equation}\label{eq66}
(\brB-H^T\brA)^T(\brB-H^T\brA)\preceq \sum_k\lambda_kS_k,
\end{equation}
for $x\in \X$ it holds
\begin{equation}\label{eq661}
\|Bx-H^TAx\|_2^2\leq\phi_{\T}(\lambda).
\end{equation}
\begin{quote}
Indeed, in the case of (\ref{eq66}) and with $x\in \X$, there exists $t\in\T$ such that $x^TS_kx\leq t_k$ for all $t$, and consequently the vector $\bar{t}$ with the entries $\bar{t}_k=x^TS_kx$ also belongs to $\T$, whence
$$
\|{Bx}-H^T{Ax}\|_2^2=\|{Bx}-H^T{Ax}\|_2^2\leq \sum_k\lambda_k{x}^TS_k{x}=\lambda^T\bar{t}\leq\phi_\T(\lambda),
$$
which combines with (\ref{Xequals}) to imply (\ref{eq661}).
\end{quote}
From (\ref{eq661}) it follows that if $H$ and $\lambda\geq0$ are linked by (\ref{eq66}), then
\bse
\Risk^2[\widehat{x}_H|\X]&=&\max_{x\in \X} \bE\left\{\|Bx-H^T(Ax+\sigma\xi)\|_2^2\right\}\\&=&
\sigma^2\Tr(H^TH)+\max_{x\in \X} \|[B-H^TA]x\|_2^2\\
&\leq&\sigma^2\Tr(H^TH)+\phi_{\T}(\lambda).\\
\ese We see that the efficiently computable convex function
$$
\widehat{R}(H)=\inf_\lambda\left\{\begin{array}{l}\sigma^2\Tr(H^TH)+\phi_{\T}(\lambda):\\~~~~~~~(\brB-H^T\brA)^T(\brB-H^T\brA)\preceq \sum_k\lambda_kS_k,\lambda\geq0\end{array}\right\}
$$
(which clearly is well defined due to compactness of $\T$ combined with $\sum_kS_k\succ0$) is an upper bound on $R(H)$.{\footnote{It is well known that when $K=1$ (i.e., $\X$ is an ellipsoid), the above bounding scheme is exact: $R(\cdot)\equiv \widehat{R}(\cdot)$. For more complicated $\X$'s, $\widehat{R}(\cdot)$ could be larger than $R(\cdot)$, although the ratio $\widehat{R}(\cdot)/R(\cdot)$ is bounded by $O(\log(K))$, see section \ref{secbypro}.}}
 Therefore the efficiently computable optimal solution  $(H_*,\lambda_*)$
to the (clearly solvable) convex optimization problem
\begin{equation}\label{eq555}
\begin{array}{rcl}
\Opt&=&\min_{H,\lambda}\left\{\sigma^2\Tr(H^TH)+\phi_{\T}(\lambda): \right.\\&&~~~~~~~~~~~\left.(\brB-H^T\brA)^T(\brB-H^T\brA)\preceq \sum_k\lambda_kS_k,\lambda\geq0\right\}
\\
&=&\min_{H,\lambda}\left\{\sigma^2\Tr(H^TH)+\phi_{\T}(\lambda):\right.\\
&&~~~~~~~~~~~\left.\left[\begin{array}{cc}
\sum_k\lambda_k S_k&\brB^T-\brA^TH\cr
\brB-H^T\brA&I_\nu\cr\end{array}\right]\succeq0,\lambda\geq0\right\}\\
\end{array}
\end{equation}
yields a linear estimate $\widehat{w}_{H_*}$ with the risk upper-bounded by $\sqrt{\Opt}$.

\subsection{Lower-bounding optimal risk and near-optimality of ${\wh{w}}_{H_*}$}\label{mainobs}
\def\brB{B}
\def\brA{A}
Let us consider the convex optimization problem
\begin{equation}\label{diezrho}
\begin{array}{l}
\Opt_*=\max\limits_{Q,t}\Big\{\varphi(Q):=\Tr\big(B[Q-QA^T(\sigma^2 I_m+AQA^T)^{-1}AQ]B^T\big),\\
\multicolumn{1}{r}{Q\succeq 0,\,t\in \cT,\,\Tr(QS_k)\leq t_k,\,1\leq k\leq K\Big\}}
\end{array}%
\end{equation}
\begin{equation}\label{starrho}
~~~~~~~~=\max\limits_{Q,t}\left\{
\begin{array}{l}\Tr(BQB^T)-\Tr(G):
\left[\begin{array}{cc}G&BQA^T\\A QB^T&\sigma^2I_m+A QA^T\end{array}\right]\succeq0,\\
~~~~~~~~~Q\succeq 0,\,t\in \T,\,\Tr(QS_k)\leq t_k,\,1\leq k\leq K
\end{array}
\right\}.
\end{equation}
Note that the function $\varphi(Q)$  has a transparent statistical interpretation. Specifically, given an $n\times n$ matrix $Q\succeq0$, consider two independent Gaussian random vectors,
$\xi\sim\N(0,I_m)$ and $\eta\sim\N(0,Q)$. We claim that
\be
\varphi(Q)=\inf_{G(\cdot):\,\bR^m\to \bR^\nu}\bE_{[\xi,\eta]}\{\|G(\sigma \xi+A\eta)-B\eta\|_2^2\}.
\ee{blueeq:gaussopt}
Indeed, by the Normal Correlation theorem (see, e.g., \cite[Theorem 13.1]{liptser1977statistics}), the optimal, in terms of expected $\|\cdot\|_2^2$-error, recovery $G_*(\cdot)$ of $B\eta$ via observation
$\sigma\xi+\brA\eta$ -- the conditional, given $\sigma \xi+A\eta$, expectation of $B\eta$ -- is linear, and the corresponding expected $\|\cdot\|_2^2$-error is exactly $\varphi(Q)$.
\par
In the sequel, we set
\[
\Q=\{Q\in\bS^n:Q\succeq0,\exists t\in\T: \Tr(QS_k)\leq t_k,\,1\leq k\leq K\}.
\]
Note that $\Q$ is a convex compact set due to $\sum_kS_k\succ0$ combined with compactness of $\T$.
\par
Observe that if $(Q,t)$ is feasible for \rf{diezrho}, then the Gaussian random vector $\eta\sim \N(0,Q)$ belongs to $\X$ ``on average'' -- it satisfies the constraints $\bE\{\eta^TS_k\eta\}=\Tr(QS_k)\leq t_k,\,k=1,...,K,$ and $t\in \T$. The lower bounding scheme we intend to implement goes back to  \cite{Pinsker1980} and heavily relies upon this fact -- it bounds from below the minimax, over $x\in \X$, risk of estimating $Bx$ by comparing this risk to the risk of optimal recovery of $B\eta$ in the Gaussian problem, where $\eta\in \X$ with ``high probability,'' as is the case when 
$Q\in \rho\Q$ with appropriate $\rho<1$. Specifically, we have the following simple
\begin{lemma}\label{lem:lower1}
 Given a positive semidefinite $n\times n$ matrix $Q$ and $\delta\in(0,1/5]$, let $\eta\sim \N(0,Q)$ and $\xi\sim \N(0,I_m)$ be independent from each other Gaussian vectors. Assume that
 $$\Prob\{\eta\notin \X\}\leq \delta\leq 1/5.$$
 Then
\be
\varphi(Q)\leq \Riskopt^2[\X]+\left[M_*+\sqrt{2} q_{1-\delta/2} {\sqrt{\Tr(BQB^T)}}\right]^2\delta,
\ee{eq:0002}
where $q_\alpha$ is the $\alpha$-quantile of the standard normal distribution: \[{1\over\sqrt{2\pi}}\int_{-\infty}^{q_\alpha}{\rm e}^{-s^2/2} ds = \alpha,\] and
\be
M_*=\sqrt{\max\limits_{Q\in\Q}\Tr(BQB^T)}.
\ee{blueM_*}
{Further}, if $Q\in \rho \Q$ for some $\rho\in(0,1]$,  then
\[
\varphi(Q)\leq \Riskopt^2[\X]+[1+\sqrt{2\rho} q_{1-\delta/2}]^2M_*^2\delta.
\]
\end{lemma}
For proof, see Section {A.1.1}.\par
The second principal component of the construction of the lower bound for $\Riskopt$ is provided by the following statement:
\begin{lemma}\label{Opto}
In the premise of this section {\rm \rf{starrho}} is a conic problem which is strictly feasible and solvable, with the conic dual problem equivalent to {\rm \rf{eq555}}. As a consequence, one has
\be
\Opt_*=\Opt.
\ee{miracle}
\end{lemma}
Let now $(Q,t)$ be an optimal solution to \rf{diezrho}, and let for $0<\rho\leq1$, $Q_\rho=\rho Q$. Note that $\varphi(Q_\rho)\geq \rho\varphi(Q)=\rho\Opt$ {(recall that $\varphi$ is concave\footnote{Concavity of $\varphi$ can be verified directly; a transparent alternative verification is to notice that (\ref{blueeq:gaussopt}) implies that
$$\begin{array}{rcl}
\varphi(Q)&=&\min_H\bE_{[\xi,\zeta]\sim\cN(0,I_{m+n})}\left\{([B-H^TA]Q^{1/2}\zeta+H^T\sigma\xi)^T([B-H^TA]Q^{1/2}\zeta+H^T\sigma\xi)\right\}\\
&=&
\min_H\left[\Tr([B-H^TA]^TQ[B-H^TA])+\sigma^2\Tr(H^TH)\right]\\
\end{array}
$$
is a minimum of affine functions of $Q$ and as such is concave.} with $\varphi(0_{m\times m})=0$)}, and
\[
\Tr(BQ_{\rho}B^T)=\rho\Tr(BQB^T)\leq  \rho M_*^2.
\]
In view of Lemma \ref{lem:lower1} as applied with $Q_\rho$ in the role of $Q$, whenever $\rho\in(0,1]$ and there exists $\delta_\rho\leq 1/5$ such that $\Prob_{\eta\sim\N(0,Q_\rho)}\{\eta\not\in \X\}\leq \delta_\rho$, we have
\begin{equation}\label{eq:00010}
\rho\Opt\leq \varphi(Q_\rho)\leq \Riskopt^2[\X]+[1+\sqrt{2\rho}\,q_{1-\delta_\rho/2}]^2M_*^2\delta_\rho.
\end{equation}
To proceed, we need an upper bound $\delta_\rho$ on the probability $\Prob_{\eta\sim\N(0,Q_\rho)}\{\eta\notin \X\}$. It is given by the following simple result.
\begin{lemma}\label{lem1} Let $S$ and $Q$ be positive semidefinite $n\times n$ matrices with $\rho:=\Tr(SQ)\leq1$, and let $\eta\sim\N(0,Q)$. Then
\begin{equation}\label{then1}
\begin{array}{rcl}\Prob\left\{\eta^TS\eta>1\right\}&\leq& \inf\limits_{0\leq \gamma<\min_i(2s_i)^{-1}}\left\{\exp\left(-\half\sum_{i=1}^{n}\ln(1-2\gamma s_i)-\gamma\right)
\right\}\\
&\leq& \e^{-{1-\rho+\rho\ln(\rho)\over 2\rho}}
\end{array}
\end{equation}
where $s_i$ are the eigenvalues of $Q^{1/2}SQ^{1/2}$.
\end{lemma}
Now we are done. Indeed, note that  the matrix $Q_\rho$ satisfies $\Tr(S_kQ_\rho)\leq \rho t_k$ for some $t\in \T$; applying Lemma \ref{lem1} and taking into account (\ref{Xequals}), we conclude that
\[
\Prob_{\eta\sim\N(0,Q_\rho)}\{\eta\notin \X\}\leq \sum_{k=1}^K\Prob\{\eta^TS_k\eta >t_k\}\leq K\exp\left\{-{1-\rho+\rho\ln(\rho)\over 2\rho}\right\},
\]
so we can set
\begin{equation}\label{deltarho}
\delta_\rho:=\min\left[K\exp\left\{-{1-\rho+\rho\ln(\rho)\over 2\rho}\right\},1\right].
\end{equation}
It is straightforward to verify that with the just defined $\delta_\rho$,  for $0<\rho<1$ it holds
\begin{equation}\label{onethird}
[1+\sqrt{2\rho}q_{1-\delta_\rho/2}]^2\delta_\rho\leq \mysix K\exp\left\{-(3\rho)^{-1}\right\}.
\end{equation}
Assuming that $\delta_\rho\leq 1/5$, the latter bound combines with \rf{eq:00010} to yield
\be
\rho\,\Opt\leq \Riskopt^2[\X]+\mysix KM_*^2\exp\left\{-(3\rho)^{-1}\right\}.
\ee{eq:rhoopt}
Let us choose
\[
\bar{\rho}^{-1}=3\ln \left(\mysix KM_*^2\over \Riskopt^2[\X]\right)
\]
so that
\[
\mysix KM_*^2\exp\left\{-(3\bar{\rho})^{-1}\right\}\leq \Riskopt^2[\X].
\]
Observe that by evident reasons $M_*^2\geq \Riskopt^2[\X]$, whence $\bar{\rho}^{-1}\geq 3\ln(\mysix K)$, which in view of (\ref{deltarho}) implies that $\delta_{\bar{\rho}}\leq 1/5$, so that (\ref{eq:rhoopt})
 is applicable to $\rho=\bar{\rho}$, thus implying that
\[
\Opt\leq {2\over \bar{\rho}}\Riskopt^2[\X]=6 \ln \left(\mysix KM_*^2\over \Riskopt^2[\X]\right)\Riskopt^2[\X].
\]
Recalling that $\sqrt{\Opt}$ upper-bounds $\Risk[\widehat{w}_{H_*}|\X]$, we have arrived at our main result:
\begin{theorem} \label{prop1} The  efficiently computable linear estimate ${\wh{w}}_{H_*}(\omega)=H_*^T\omega$ yielded by an optimal solution to the optimization problem {\rm (\ref{eq555})} is nearly optimal in terms of its risk:
\begin{equation}\label{nearoptimal}
\Risk[{\wh{w}}_{H_*}|{\X}]\leq\sqrt{\Opt}\leq \sqrt{6 \ln\left({\mysix M_*^2K\over \tRiskopt^2[{\X}]}\right)}\Riskopt[{\X}]
\end{equation}
with $M_*$ given by {\rm (\ref{blueM_*})}.
\end{theorem}
\subsection{Discussion}
The result of Theorem \ref{prop1} merits few comments.
\paragraph{1. Simplifying {expression for} nonoptimality factor} Relation (\ref{nearoptimal})  states  that when ${\X}$ is an {ellitope} (\ref{setX}), the risk $\sqrt{\Opt}$ of the efficiently computable linear estimate yielded by (\ref{eq555}) is just by a logarithmic in ${{M_*^2}K\over
\tRiskopt^2[{\X}]}$ factor worse than the optimal risk $\Riskopt[{\X}]$. A {minor} shortcoming of (\ref{nearoptimal}) is that the ``nonoptimality factor'' is expressed in terms of unknown to us optimal risk. This can be easily cured. For example, setting
\def\myfourteen{{\cbr 17}}
$$
\bar{\rho}^{-1}=6  \ln \left(17 KM_*^2\over \Opt\right),
$$
it is immediately seen that
\[{\bar{\rho}\over 2}\Opt\geq 6 KM_*^2\exp\{-(3\bar{\rho})^{-1}\},
\]
and $\delta_{\bar{\rho}}$ as given by (\ref{deltarho}) with $\rho=\bar{\rho}$ is $\leq 1/5$, implying by (\ref{eq:rhoopt}) that
${1\over 2}\bar{\rho}\Opt\leq \Riskopt^2[\X],$ whence
\be
\Riskopt^2[\X]\geq \left[12\ln \left(17 KM_*^2\over \Opt\right)\right]^{-1} \Opt.
\ee{eq:optopt}
Note  that all the quantities in the right hand side of \rf{eq:optopt} are efficiently computable given the problem data, and that $\sqrt{\Opt}$ is an upper bound on $\Risk[\widehat{w}_{H_*}|\X]$.

Furthermore, if a simple though less precise expression of the factor in terms of this data is required, it can be obtained as follows.
Recall that two points $x=x_+$ and $x=-x_+$ of $\X$ can be distinguished through the observation $Ax+\sigma\xi$ with maximal probability of error $0<\alpha<1$ only if $\|Ax\|_2
   \geq c_\alpha\sigma$, $c_\alpha>0$;\footnote{In fact, one can choose $c_\alpha=q_{1-\alpha}$, the $1-\alpha$-quantile of the standard normal distribution.} by the standard argument one conclude that
    the risk of estimation of $Bx$ satisfies, for some absolute constant $c>0$:
\be
\Riskopt^2[\X]\geq \max\left\{\|Bx\|_2:\,\|Ax\|_2\leq c\sigma,\,x\in \X\right\}.
\ee{eq:slower1}
Now let $B=I$, and consider two typical for the traditional non-parametric statistics types of ${\X}$:
\begin{itemize}
\item ${\X}$ is the ellipsoid $\{x\in\bR^n: \sum_ia_i^2x_i^2\leq 1\}$ with $0<a_1\leq a_2\leq...\leq a_n$ (for properly selected $a_i$ this set models the restriction onto a regular $n$-point grid of functions from a Sobolev ball).
Here $K=1$, $\T=[0,1]$, $S_1=\Diag\{a_1^2,...,a_n^2\}$. When choosing $x=t\e_1$, where $e_1$ is the first basic orth and $t\in ]0,1]$, using \rf{eq:slower1} we get
$\Riskopt[\X]\geq \min\left[1/a_1,\, c\sigma/\|[A]_1\|_2\right]$ where $[A]_1$ is the first column of $A$. On the other hand, we have $M_*^2=a_1^{-2}$, and the simplified risk bound reads
 \[
    \Risk[{\wh{w}}_{H_*}|{\X}]\leq O(1)
    \sqrt{\ln\left(1+{\|[A]_1\|_2\over \sigma a_1}\right)}\Riskopt[{\X}].
\]\item ${\X}$ is the box $\{x\in\bR^n: a_i|x_i|\leq 1,\,1\leq i\leq n\}$, where, as above, $0<a_1\leq a_2\leq...\leq a_n$. Here $K=n$, $\T=[0,1]^n$, $x^TS_kx=a_k^2x_k^2$, resulting in
$M_*^2=\sum_{i} a_i^{-2}\leq n a_1^{-2}$.
The same bound
$\Riskopt[\X]\geq \min\left[1/a_1,\, c\sigma/\|[A]_1\|_2\right]$ holds in this case and, consequently,
\[
\Risk[{\wh{w}}_{H_*}|{\X}]\leq O(1)\sqrt{\ln n+\ln \left(1+{\|[A]_1\|_2\over \sigma a_1}\right)}\Riskopt[{\X}].
\]
    \end{itemize}
Now let $B$ be a general-type matrix, and assume for the sake of simplicity that $B$ has trivial kernel.
We associate with the data the following quantities:
\begin{itemize}
\item size of $\T$, $T=\max_{t\in \T}\sum_{k}t_k$, and $\varkappa$ -- the minimal eigenvalue of
$\sum_k S_k$.
Note that for any $x\in \X$, $\sum_{k}x^TS_kx\leq T$, thus the radius
$r(\X)=\max_{x\in \X}\|x\|_2$ of $\X$ satisfies $r(\X)\leq \sqrt{T/\kappa}$;
\item $\ell_1/\ell_\infty$-{\sl condition number of $\T$}
\[
\cond(\T)=\sqrt{{T\over \max\limits_{t\in \T}\min\limits_{k\leq K}t_k}}
=\sqrt{{\max\limits_{t\in\T}\sum_k t_k\over \max\limits_{t\in \T}\min\limits_{k\leq K}t_k}};
\]
by our assumptions, $\T$ intersects the interior of $\bR^K_+$ and thus $\sqrt{K}\leq \cond(\T)<\infty$;
\item {\sl condition number of $B$:}
$\cond(B)={\sigma_{\max}(B)\over\sigma_{\min}(B)}$,
where $\sigma_{\max}(B)$ and $\sigma_{\min}(B)$ are, respectively,  the largest and the smallest singular values of $B$.
\end{itemize}
\begin{corollary}\label{simplefactor}
In the situation of this section
\begin{equation}\label{eq:simplef}
\Risk[{\wh{w}}_{H_*}|{\X}]\leq O(1)\sqrt{\ln\left(K{\cond}^2(B)\left[{\cond}^2(\T)+{\|A\|^2T\over\sigma^2\varkappa}\right]\right)}
\Riskopt[{\X}];
\end{equation}
here and in what follows, $O(1)$ stands for a properly selected positive absolute constant.
\end{corollary}
It is worth to note that, surprisingly, the logarithmic factor in \rf{eq:simplef} does not depend of the structure of singular spectrum of $A$, the entity which, as far as the role of $A$ is concerned, is primarily responsible for $\Riskopt[{\X}]$.
\paragraph{2. Relaxing the symmetry requirement} Sets ${\X}$ of the form (\ref{setX}) -- we called them {ellitopes} -- are symmetric w.r.t. the origin convex compacts of special structure. This structure is rather flexible, but the symmetry is ``built in.'' We are about to demonstrate that, to some extent, the symmetry requirement can be relaxed. Specifically, assume instead of (\ref{setX}) that for some $\alpha\geq1$ it holds
\[
\underbrace{\big\{x\in\bR^n:\exists (y\in\bR^{\brn},t\in\T): x=Py\ \&\  y^TS_ky\leq t_k,\,1\leq k\leq K\big\}}_{\underline{{\X}}
}\subset {\X}\subset \alpha\underline{{\X}},
\]
with $S_k$ and $\T$ possessing the properties postulated in section \ref{sitgoal}. Let $\Opt$ and $H_*$ be the optimal value and  optimal solution of the optimization problem (\ref{eq555}) associated with the data $S_1,\, ...,S_K,\,\T$ and matrices $\bar{A}=AP$, $\bar{B}=BP$ in the role of $A$, $B$, respectively. It is immediately seen that the risk $\Risk[{\wh{w}}_{H_*}|{\X}]$ of the linear estimate ${\wh{w}}_{H_*}(\omega)$ is at most $\alpha\sqrt{\Opt}$. On the other hand,  we have $\Riskopt[\underline{{\X}}]\leq\Riskopt[{\X}]$, and by Theorem \ref{prop1} also
$\sqrt{\Opt}\leq \sqrt{6\ln\left({\mysix M_*^2K\over \tRiskopt^2[\underline{{\X}}]}\right)}\Riskopt[\underline{{\X}}]$. Taken together, these relations imply that
\begin{equation}\label{alpha}
\Risk[{\wh{w}}_{H^*}|{\X}]\leq \alpha\sqrt{6 \ln\left({{\mysix M_*^2}K\alpha\over \tRiskopt^2[{\X}]}\right)}\Riskopt[{\X}].
\end{equation}
{In other words,} as far as the ``level of nonoptimality'' of efficiently computable linear estimates is concerned, signal sets ${\X}$ which can be approximated by {ellitopes} within a factor $\alpha$ {\sl of order of 1} are nearly as good as the {ellitopes}. To give an example: it is known that whenever the intersection ${\X}$ of $K$ elliptic cylinders $\{x:(x-c_k)^TS_k(x-c_k)\leq1\}$, $S_k\succeq 0$, concentric or not, is bounded and has a nonempty interior, ${\X}$ can be approximated by an ellipsoid within the factor $\alpha=K+2\sqrt{K}$ \footnote{specifically, setting $F(x)=-\sum_{k=1}^K\ln(1-(x-c_k)^TS_k(x-c_k)):\inter {\X}\to \bR$ and denoting by $\bar{x}$ the {\sl analytic center} $\hbox{\rm argmin}_{x\in\inter {\X}}F(x)$, one has
  $$
  \{x:(x-\bar{x})^TF''(\bar{x})(x-\bar{x})\leq1\}\subset {\X}\subset \{x:(x-\bar{x})^TF''(\bar{x})(x-\bar{x})\leq [K+2\sqrt{K}]^2\}.
  $$
  }. Assuming w.l.o.g. that the approximating ellipsoid is centered at the origin, the level of nonoptimality of a linear estimate is bounded by (\ref{alpha}) with $O(1)K$ in the role of $\alpha$. Note that  bound (\ref{alpha}) rapidly deteriorates when $\alpha$ grows,
  and this phenomenon to some extent ``reflects the reality.'' For example, a {perfect} simplex ${\X}$ inscribed into the unit sphere in $\bR^n$ is in-between two centered at the origin Euclidean balls with the ratio of radii equal to $n$ (i.e. $\alpha=n$). It is immediately seen that with $A=B=I$, in the range $\sigma\leq n\sigma^2\leq1$ of values of $n$ and $\sigma$, we have
  $$
  \Riskopt[{\X}]\approx \sqrt{\sigma},\,\,\Riskopt[{\wh{w}}_{H_*}|{\X}]=O(1) \sqrt{n}\sigma,$$
  with $\approx$ meaning ``up to logarithmic in $n/\sigma$ factor.''
   In other words, for large $n\sigma$ linear estimates indeed are significantly (albeit not to the full extent of (\ref{alpha})) outperformed by nonlinear ones.
\par Another ``bad for linear estimates'' situation suggested by (\ref{nearoptimal}) is that where the description (\ref{setX}) of ${\X}$, albeit possible, requires {a huge value} of $K$. Here again (\ref{nearoptimal}) reflects to some extent the reality: when ${\X}$ is the unit $\ell_1$ ball in $\bR^n$, (\ref{setX}) takes place with $K=2^{n-1}$; consequently, the factor at $\Riskopt[{\X}]$ in the right hand side of (\ref{nearoptimal}) becomes at least $\sqrt{n}$. On the other hand, in the range $\sigma\leq n\sigma^2\leq1$ of values of $n$, $\sigma$, and with $A=B=I$, the risks $\Riskopt[{\X}]$, $\Riskopt[{\wh{w}}_{H_*}|{\X}]$ are basically the same as
in the case of ${\X}$ being the perfect simplex inscribed into the unit sphere in $\bR^n$, and linear estimates indeed are ``heavily non-optimal'' when $n\sigma$ is large.

\subsection{Numerical illustration}\label{sec:num1}
Observe that inequality \rf{eq:00010}
taken together with an efficiently computable
upper bound $\delta_\rho$ for the probability that $\eta\notin \X$ for $\eta \sim\N(0,Q_\rho)$ yields a single-parametric family of lower bounds on $\Riskopt[{\X}]$:
\[
\Riskopt^2[X]\geq\rho\,\Opt-[1+\sqrt{2\rho}\,q_{1-\delta_\rho/2}]^2M_*^2\delta_\rho.
\]
We can compute the right hand side for several values of $\rho$, take the largest of the resulting lower bounds on $\Riskopt[{\X}]$ and compare the result with the risk $\sqrt{\Opt}$ of
the efficiently computable linear estimate yielded by the optimal solution to (\ref{eq555}). In this way, we hopefully will end up with less pessimistic assessment
of the level of non-optimality of linear estimates than the one yielded by \rf{nearoptimal}. On the other hand, better lower bounds can be computed using directly the inequality \rf{eq:0002} of Lemma \ref{lem:lower1} along with an efficiently computable  approximation of the constraint
$\Prob\{\eta\notin \X\}\leq \delta$ on the distribution $\N(0,Q)$ of $\eta$. Indeed, given $0<\delta\leq1/5$, suppose that ${\Q}_\delta$ is a convex subset of the positive semidefinite cone such that for any $Q\in{\Q}_\delta$ and $\eta\sim \N(0,Q)$ one has
$\Prob\{\eta\notin \X\}\leq \delta$. Then, according to \rf{eq:0002},  the quantity
\be
\Opt_\delta-[M_*+\sqrt{2} q_{1-\delta/2} \|BQ_\delta^{1/2}\|_{2}]^2\delta,
\ee{eq:lrb}
where
\[
\Opt_\delta=\max_{Q\in {\Q}_\delta} \varphi(Q)
\]
and $Q_\delta$ is the corresponding optimal solution, is a lower bound on $\Riskopt[{\X}]$.
\par
\begin{figure}[h]
\begin{tabular}{cc}
\includegraphics[scale=0.4]{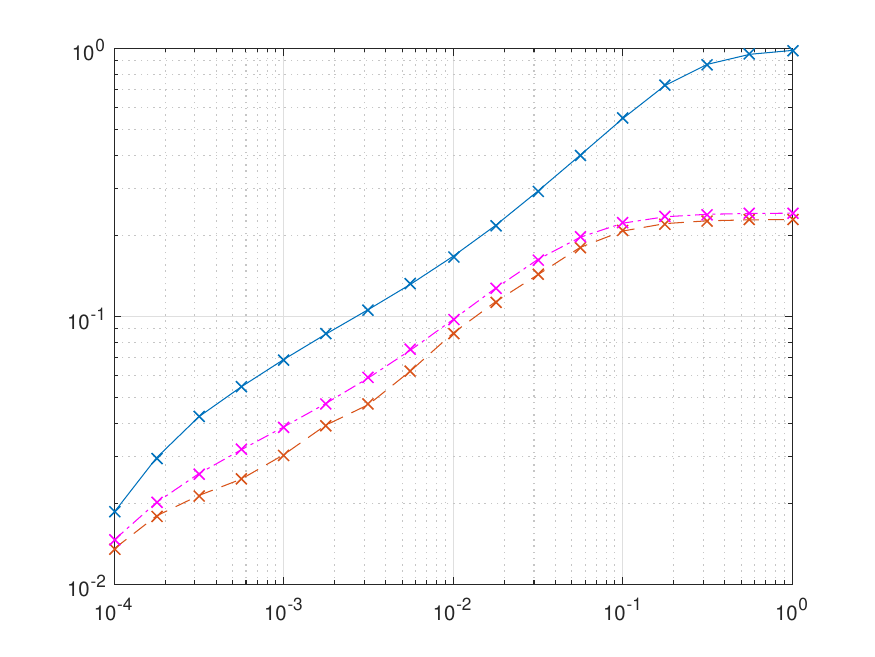}&\includegraphics[scale=0.4]{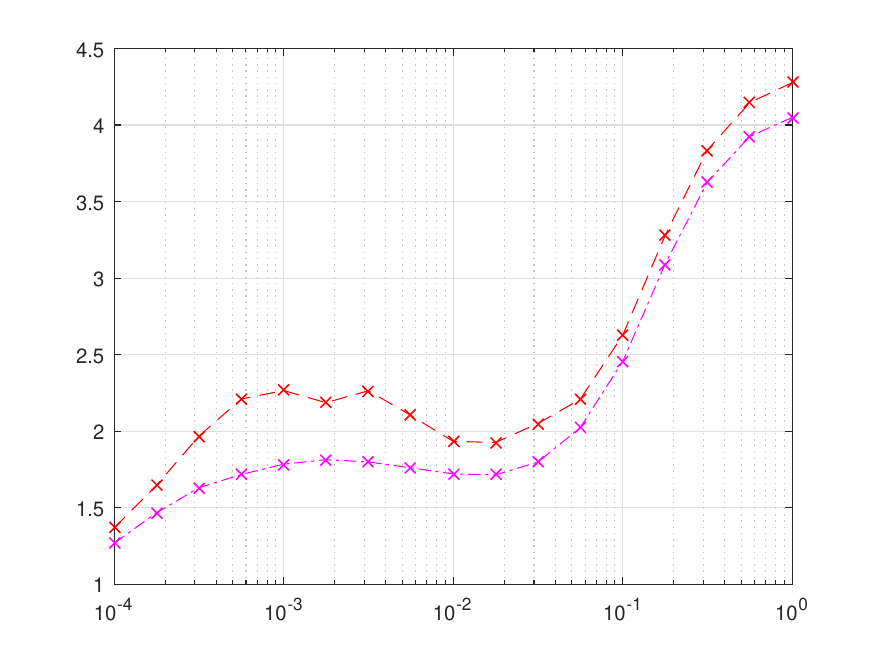}
\end{tabular}
\caption{\small Recovery on ellipsoids: risk bounds as functions of the noise level $\sigma$, dimension $n=32$. Left plot: upper bound of the risk of linear recovery (solid blue line); red dash line and magenta dash-dot line -- lower bounds utilizing two implementations of the bounding scheme.  Right plot: suboptimality ratios.}
\label{fig:e1-1}
\end{figure}
\begin{figure}[h]

\begin{tabular}{cc}
\includegraphics[scale=0.4]{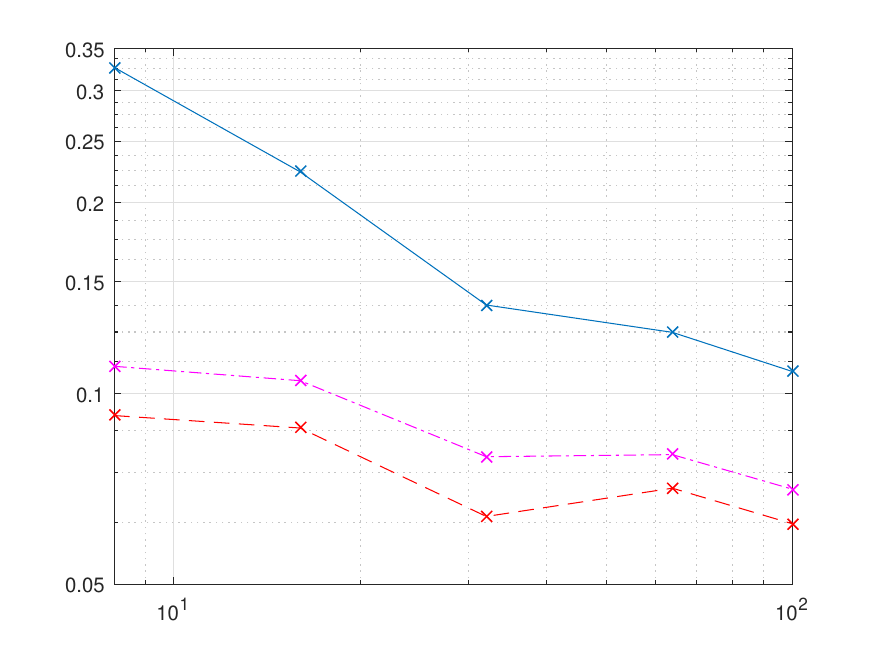}&\includegraphics[scale=0.4]{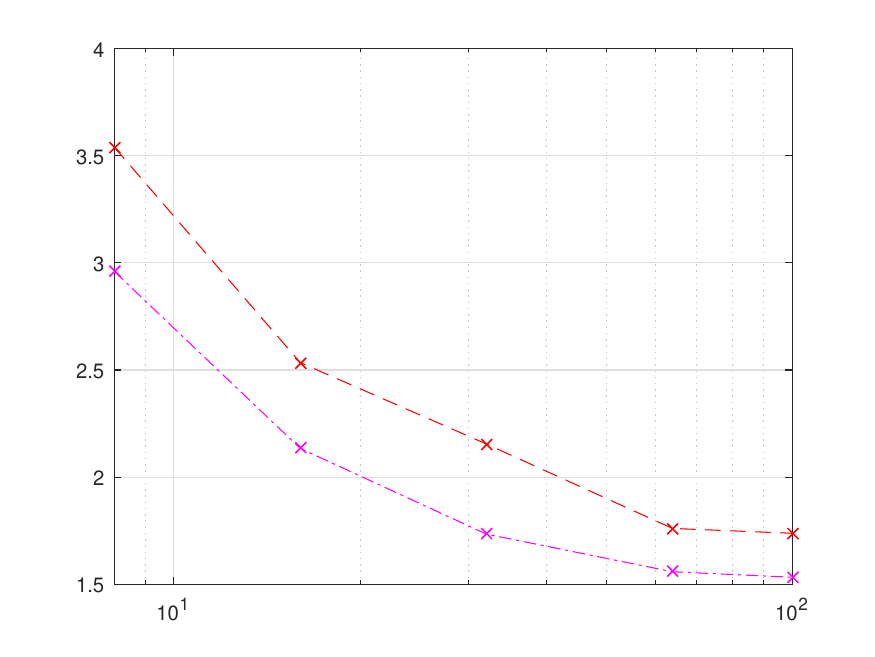}
\end{tabular}
\caption{\small Recovery on ellipsoids: risk bounds as functions of problem dimension $n$, noise level $\sigma=0.01$. Left plot: upper bound of the risk of linear recovery (solid blue line), red dash line and magenta dash-dot line -- lower bounds on $\Riskopt$ utilizing two implementations of the bounding scheme. Right plot:  suboptimality ratios.}
\label{fig:e1-end}
\end{figure}

We have conducted two experiments aimed to compare the sub-optimality factors obtained numerically with their theoretical counterparts. In both experiments $B$ and $P$ are set to be $n\times n$ identity matrices, and $n\times n$ sensing matrix $A$ is a randomly rotated matrix with singular values $\lambda_j$, $1\leq j\leq n$, forming a geometric progression, with $\lambda_1=1$ and $\lambda_{n}=0.01$. In the first experiment the signal set ${\X}_1$ is an ellipsoid:
\[
{\X}_1=\{x\in\bR^{n}: \sum_{j=1}^{n}j^2x_j^2\leq1\},
\]
that is, $K=1$, $S_1=\sum_{j=1}^{n}j^2e_je_j^T$ ($e_j$ are basic orths), and $\T=[0,1]$.
With two natural implementations of the outlined bounding scheme (for the sake of completeness, the details of the lower bound computation are provided in {Section C},
we arrived at simulation results presented on Figures \ref{fig:e1-1} and \ref{fig:e1-end}. {It is worth to mention that the theoretical estimation of the ``suboptimality factor'' computed according to \rf{eq:optopt} varies in the interval $[31.6,\,73.7]$ in this experiment.}

\begin{figure}[ht]

\begin{tabular}{cc}
\includegraphics[scale=0.4]{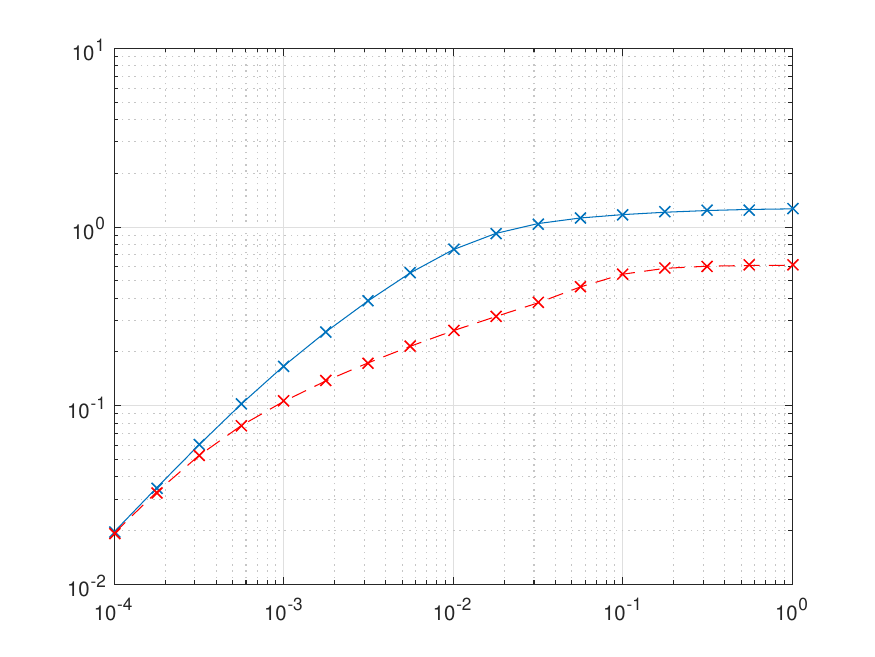}&\includegraphics[scale=0.4]{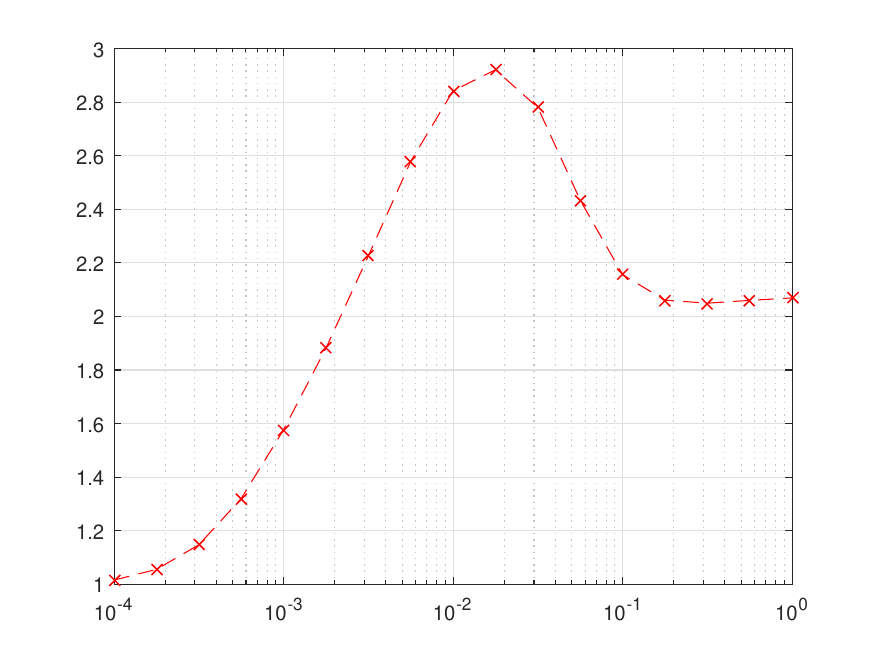}
\end{tabular}
\caption{\small Recovery on a box: risk bounds as functions of the noise level $\sigma$, dimension $n=32$. Left plot: upper bound of the risk of linear recovery (solid blue line) and lower risk bound (red dash line). Right plot: suboptimality ratios.}
\label{fig:e2-1}
\end{figure}

\begin{figure}[ht]

\begin{tabular}{cc}
\includegraphics[scale=0.4]{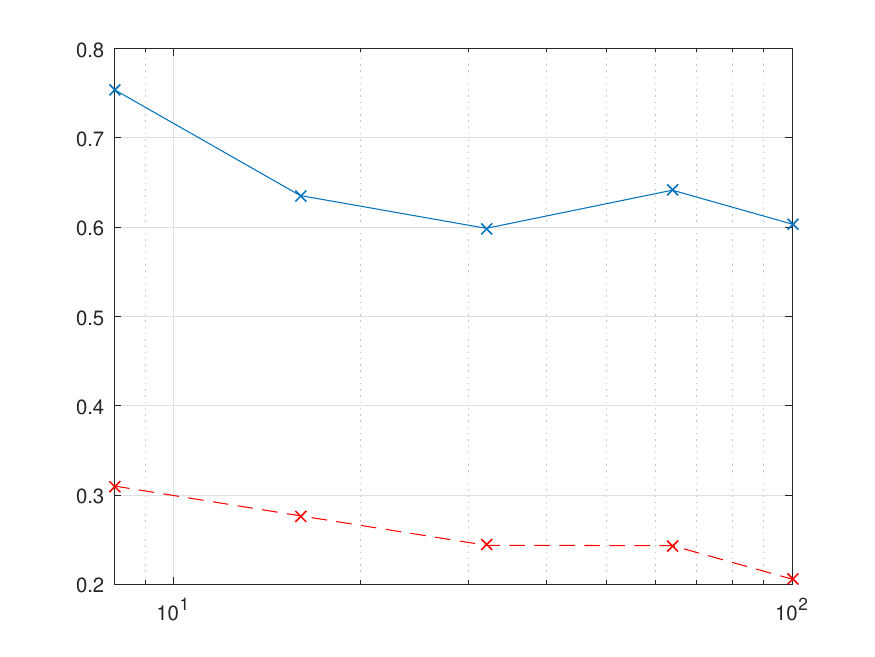}&\includegraphics[scale=0.4]{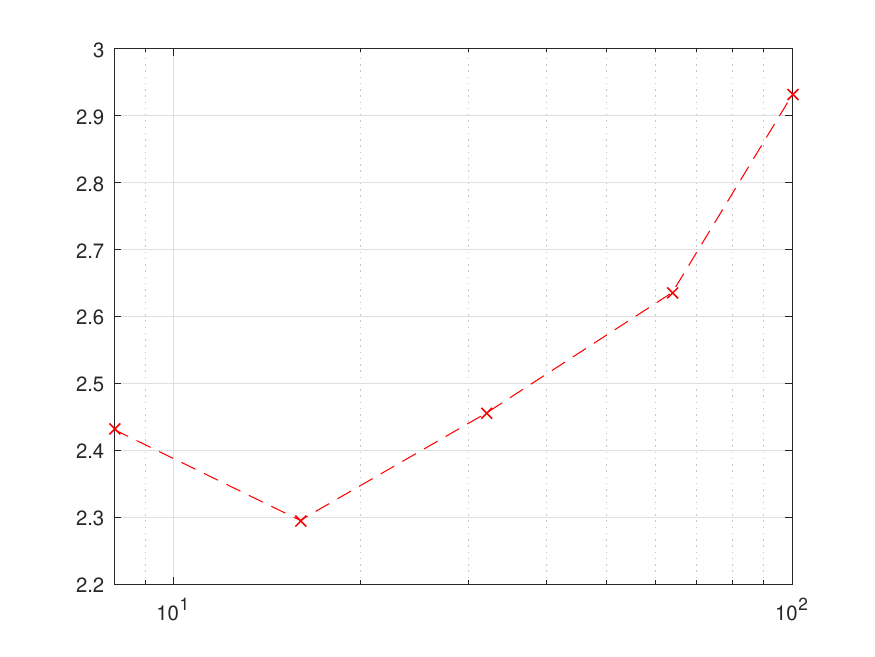}
\end{tabular}
\caption{\small Recovery on a box: risk bounds as functions of problem dimension $n$, noise level $\sigma=0.01$. Left plot: upper bound of the risk of linear recovery (solid blue line) and lower  bound on $\Riskopt$ (red dash line). Right plot: suboptimality ratio.}
\label{fig:e2-end}
\end{figure}

In the second experiment, the signal set $\X$ is the box circumscribed around the above ellipsoid:
\[
{\X}=\{x\in\bR^n: j|x_j|\leq1,\,1\leq j\leq n\}\eqno{[K=n,S_k=k^2e_ke_k^T,k=1,..., K,\T=[0,1]^K]}.
\]
In this case only one implementation of the bounding scheme is used. The simulation results of the second experiment are given on Figures \ref{fig:e2-1} and \ref{fig:e2-end}. {In this experiment also, the theoretical estimation of the non-optimality of the linear estimate is very conservative -- for different values of parameters the factor in the bound \rf{eq:optopt} varies between 73.2 and 115.4.}
\section{Extensions}\label{sectextensions}
\subsection{Estimation in relative scale}\label{estrelscale}
In this section we consider the setting as follows.
Assume that, same as in section \ref{sect1}, we are given a $\nu\times n$ matrix $B$, and a noisy observation
\[\omega=Ax+\sigma\xi,\;\;\xi\sim\N(0,I_m),
\]
of a signal $x\in \X$
with known $m\times n$ matrix $A$ and $\sigma>0$, and we aim to recover $w=Bx$.
We are given a positive semidefinite symmetric $n\times n$ matrix $S$, and we quantify the quality of a candidate estimate ${\wh{w}}(\cdot)$ by its {\sl $S$-risk} -- the quantity
\begin{equation}\label{S-riskdef}
\RiskS[{\wh{w}}|\X]=\inf\left\{\sqrt{\tau}:\;
\bE\{\|{\wh{w}}(Ax+\sigma\xi)-Bx\|_2^2\}\leq \tau(1+x^TSx)\;\forall x\in \X\right\}.
\end{equation}
The $S$-risk can be seen as risk with respect to the scale given by the ``regularity parameter'' $x^TSx$ of the unknown signal $x$. In particular, when $S=B^TB$, squared $S$-risk can be thought of as {\sl relative} risk -- the worst, over $x\in \bR^n$, expected $\|\cdot\|_2^2$-error of recovering $Bx$ scaled by $\|Bx\|_2^2$; when $S=0$, we arrive at the usual risk $\Risk[{\wh{w}}|{\X}]$.
\par
Same as in section \ref{sect1}, we assume w.l.o.g. that $\X$ is an ellitope given by (\ref{Xequals}) \footnote{To reduce the general case (\ref{setX}) to this one with $P=I$ it suffices to ``lift'' $A$, $B$, $S$ to the $y$-space according to $A\mapsto\bar{A}=AP$, $B\mapsto\bar{B}=BP$, $S\mapsto\bar{S}=P^TSP$ and then replace $\X$ with the set $\Y=\{y\in\bR^{\bar{n}}: \exists t\in\T: y^TS_ky\leq t_k,1\leq k\leq K\}$.}. Besides this, we assume that $B\neq0$ -- otherwise the estimation problem is trivial.\par
We are about to prove that in the situation in question, efficiently computable linear estimate is near-optimal.
\subsubsection{Building linear estimate}
Given a linear estimate ${\wh{w}}_H(\omega)=H^T\omega$ and $\tau\geq0$, let $\lambda\geq0$ be such that $[B-H^TA]^T[B-H^TA]\preceq \sum_k \lambda_kS_k+\tau S$, see (\ref{setX}), implying that for all $x\in \X$, there exists $t=t_x\in\T$ such that
\[\begin{array}{rcl}
    \bE_\xi\{\|{\wh{w}}_H(Ax)-Bx\|_2^2\}&\leq &x^T[{\sum}_k\lambda_kS_k+\tau S]x+\sigma^2\Tr(H^TH)\\&\leq& \sum_{k}t_k\lambda_k+x^TSx+\sigma^2\Tr(H^TH),\end{array}
\]
so that  for all $x\in \X$
\[
\bE_\xi\{\|{\wh{w}}_H(Ax+\sigma\xi)-Bx\|_2^2\}\leq \phi_{\T}(\lambda)+\tau x^TSx+\sigma^2\Tr(H^TH),
\]
where $\phi_{\T}$ is the support function of $\T$.
As a result, whenever $H$, $\tau\geq0$ and $\lambda\geq0$ are such that
\[
\sigma^2\Tr(H^TH)+\phi_{\T}(\lambda)\leq\tau,\;\;
(H^TA-B)^T(H^TA-B)
\preceq \sum_k\lambda_kS_k+\tau S,
\]
we have
$$
\RiskS[{\wh{w}}_H|\X]\leq\sqrt{\tau}.
$$
We arrive at the convex problem
\begin{equation}\label{designproblem}
\Opt=\min_{\tau,H,\lambda}\left\{\begin{array}{rcl}\tau:\left[\begin{array}{cc}\sum_k\lambda_k S_k+\tau S&B^T-A^TH\cr
B-H^TA&I_\nu\cr\end{array}\right]\succeq0,\\
~~~~~~~~\sigma^2\Tr(H^TH)+\phi_{\T}(\lambda)\leq\tau,\;\lambda\geq0\end{array}\right\}.
\end{equation}
The $H$-component $H_*$ of an optimal solution to this problem yields linear estimate ${\wh{w}}_{H_*}(\omega)=H_*^T\omega$ with $S$-risk $\leq \sqrt{\Opt}$.

\subsubsection{Lower-bounding the optimal $S$-risk and near-optimality of ${\wh{w}}_{H_*}$}
Consider the problem
\begin{equation}\label{eq:diez}
\Opt_*=\max_{W,G,s,v}\left\{\Tr(BWB^T)-\Tr(G):\begin{array}{l}
\left[\begin{array}{cc}G&BWA^T\cr A WB^T&\sigma^2sI_m+A WA^T\cr\end{array}\right]\succeq0,\\
W\succeq0,\,\Tr(WS_k)\leq v_k,\,1\leq k\leq K,\\
\Tr(WS)+s\leq1,\,[v;s]\in\bT\\
\end{array}\right\}
\end{equation}
where {
\begin{equation}\label{conebT}
\bT=\cl \{[t;\tau]\in\bR^K\times\bR:\tau>0,\tau^{-1}t\in\T\}\subset\bR^{K+1}_+
\end{equation}
is a closed and pointed convex cone in $\bR^{K+1}$ with a nonempty interior.
We have the following counterpart of Lemma \ref{Opto} for the present setting.
\begin{lemma}\label{OptSo}
Problem {\rm\rf{eq:diez}} is strictly feasible and solvable. Furthermore, if $(W,G,[v;s])$ is an optimal solution to {\rm \rf{eq:diez}}, then $s>0$, and
\be
\Opt=\Opt_*=
\Tr\big(B[W-WA^T(\sigma^2s I_m +A WA^T)^{-1}A W]B^T\big).
\ee{equal}
\end{lemma}}
Now let $W,v$ and $s$ stem from an optimal solution to \rf{eq:diez}. Then, as we have seen, $s>0$, and we can set $t=v/s$, so that $t\in\T$. Let also $\rho\in(0,1]$, and let us put  $Q_\rho=\rho W/s$ and $\eta\sim \N(0,Q_\rho)$. We have $S^{-1}W\succeq0$ and  $\Tr(s^{-1}W S_k)\leq t_k$, $k\leq K$, so that $s^{-1}W\in \Q$ and therefore $Q_\rho\in\rho\Q$. Hence, same as in the case of the usual risk, by Lemma \ref{lem1},
\begin{equation}\label{deltabarrhoagain}
\Prob\{\eta\not\in \X\}\leq \delta_\rho:=\min\left[K\exp\left\{-{1-\rho+\rho\ln(\rho)\over 2\rho}\right\},1\right].
\end{equation}
We also have the following analog of Lemma \ref{lem:lower1}:
\begin{lemma}\label{lem:lower2}
Given $\rho\in(0,1]$, $Q\in \rho \Q$ and $\delta\leq 1/5$, let $\eta\sim \N(0,Q)$ and $\xi\sim \N(0,I_m)$ be independent from each other Gaussian vectors. Assume that
 $$\Prob\{\eta\notin \X\}\leq \delta.$$
 Then
\be
\varphi(Q)\leq \RiskoptS^2[\X](1+\Tr(QS))+[1+\sqrt{2\rho} q_{1-\delta/2}]^2M_*^2\delta,
\ee{eq:000222}
where $M_*$ is given by {\rm (\ref{blueM_*})}, $q_\alpha$, same as in Lemma \ref{lem:lower1},  is the $\alpha$-quantile of the standard normal distribution, and
\[
\RiskoptS[\X]=\inf\limits_{\wh{w}(\cdot)}\RiskS[\wh{w}|\X].
\]
is the minimax $S$-risk associated with $\X$.
\end{lemma}
For proof, see Section {A.1.1}.\par
Now note that
\[\begin{array}{rcl}
\varphi(Q_\rho)&=& \Tr\big(B[Q_\rho-Q_\rho A^T(\sigma^2I_m+AQ_\rho A^T)^{-1}AQ_\rho]B^T\big)\\
&=&{\rho\over s}
\Tr\big(B[W-\rho W A^T(s\sigma^2I_m+\rho AWA^T)^{-1}AW]B^T\big)\\
&\geq& {\rho\over s}\Opt_*={\rho\over s}\Opt
\end{array}\]
(we have used (\ref{equal}) and the positivity of $s$).
Thus, when applying Lemma \ref{lem:lower2} with $Q_\rho$ and $\delta_\rho$ in the role of $Q$ and $\delta$,
 we obtain for all $0<\rho\leq 1$ such that  $\delta_\rho\leq 1/5$:
\begin{equation}\label{bebebe}
\begin{array}{rcl}
{\rho\over s}\Opt&\leq& \RiskoptS^2[\X]\left(1+\Tr(Q_\rho S)\right)+\left[1+\sqrt{2\rho}q_{1-\delta_\rho/2}\right]^2M_*^2\delta_\rho\\
&=&\RiskoptS^2[\X]\left(1+{\rho\over s}\Tr(W S)\right)+\left[1+\sqrt{\rho}q_{1-\delta_\rho/2}\right]^2M_*^2\delta_\rho.
\end{array}
\end{equation}
Similarly to section \ref{mainobs}, setting
\[
\bar{\rho}^{-1}=3\ln \left(\mysix KM_*^2\over \RiskoptS^2[\X]\right)
\]
we ensure that
\[
\mysix KM_*^2\exp\left\{-(3\bar{\rho})^{-1}\right\}\leq \RiskoptS^2[\X].
\]
Now, same as in the case of usual risk, we clearly have $M_*^2\geq \RiskoptS^2[\X]$, whence $\delta_{\bar{\rho}}\leq \exp\{-{1\over 3\bar{\rho}}\}\leq 1/5$, see (\ref{deltabarrhoagain}), so that (\ref{bebebe})
 is applicable with $\rho=\bar{\rho}$, thus implying that
\[
{\bar{\rho}\over s}\Opt\leq \RiskoptS^2[\X]\left(1+{\bar{\rho}\over s}\Tr(W S)\right)+\mysix KM_*^2\exp\{-{1\over 3\bar{\rho}}\},
 \]
 and
 \bse
 \bar{\rho}\Opt&\leq& \RiskoptS^2[\X]\left(s+\bar{\rho}\Tr(W S)\right)+\mysix sKM_*^2\exp\{-{1\over 3\bar{\rho}}\}\\&\leq&
 \RiskoptS^2[\X]+\mysix KM_*^2\exp\{-{1\over 3\bar{\rho}}\}=2\RiskoptS^2[\X]\\
\ese
 (note that $s+\bar{\rho}\Tr(WS)\leq s+\Tr(WS)\leq 1$ by constraints in (\ref{eq:diez})).
Recalling that $\sqrt{\Opt}$ upper-bounds $\RiskS[\widehat{w}_{H_*}|\X]$, we arrive at the following
\begin{proposition}\label{prop2} The  efficiently computable linear estimate ${\wh{w}}_{H_*}(\omega)=H_*^T\omega$ yielded by an optimal solution to the optimization problem in {\rm (\ref{designproblem})} is nearly optimal in terms of $S$-risk:
\[
\begin{array}{c}
\RiskS[{\wh{w}}_{H_*}|\X]\leq \sqrt{6 \ln\left({{\mysix }K{M^2_*}\over \tRiskoptS^2[X]}\right)}\RiskoptS[X],
\end{array}
\]
where $M_*$ is given by {\rm\rf{blueM_*}}.
\end{proposition}
\subsubsection{The case of $\X=\bR^n$}
The problem of minimizing the worst-case, over $x\in \X$, $S$-risk over linear/all possible estimates makes sense for unbounded  $\X$'s as well as for bounded ones. We intend to consider the case where $\X=\bR^n$ and to show that in this case an efficiently computable linear estimate is {\sl exactly optimal.}
\par
Similar to (\ref{designproblem}), the problem of building the best, in terms of its worst-case over $x\in\bR^n$ $S$-risk, linear estimate reads
\begin{equation}\label{designproblementirespace}
\Opt=\min_{\tau,H}\left\{\tau:\left[\begin{array}{cc}\tau S&B^T-A^TH\cr
B-H^TA&I_\nu\cr\end{array}\right]\succeq0,\;\sigma^2\Tr(H^TH)\leq\tau\right\};
\end{equation}
a feasible solution $(\tau,H)$ to this problem produces an estimate $\widehat{w}_H(\omega)=H^T\omega$ with $\RiskS[\widehat{w}_H|\bR^n]\leq\sqrt{\tau}$. We are about to demonstrate that
\begin{proposition}\label{prop6666}
Assuming problem {\rm (\ref{designproblementirespace})} feasible, the problem is solvable, and its optimal solution $(\Opt,H_*)$ induces linear estimate $\widehat{w}_{H_*}$ which is minimax optimal:
\begin{equation}\label{minmaxoptimality}
\RiskS[\widehat{w}_{H_*}|\bR^n]=\sqrt{\Opt}=\inf_{\widehat{w}(\cdot)}\RiskS[\widehat{w}(\cdot)|\bR^n].
\end{equation}
\end{proposition}
It may be interesting to compare the optimal $S$-risk $\RiskS[\widehat{w}_{H_*}|\bR^n]=\sqrt{\Opt}$ to the maximal risk $\Risk[\wh{w}_{H^*}|\X_S]$
of the optimal linear estimation of $Bx$ over the ellipsoid $\X_S=\{x\in \bR^n:\,x^TSx\leq 1\}$, so that $H^*$ is the optimal solution to \rf{eq555} with $K=1$, $S_1=S$ and $\T=[0,1]$;  note that in this case the optimal value in (\ref{eq555}) is exactly $\Risk[\wh{w}_{H^*}|\X_S]$, and not just an upper bound on this risk.  When comparing \rf{eq555} with  \rf{designproblementirespace} one can easily see that both risks are equivalent up to a factor $\sqrt{2}$:
\[
\RiskS[\widehat{w}_{H_*}|\bR^n]\leq \Risk[\widehat{w}_{H^*}|\X_S]\leq \sqrt{2}\RiskS[\widehat{w}_{H_*}|\bR^n].
\]
Note also that by the definition of $S$-risk, we have \[\Risk[\wh{w}_{H_*}|\X_S]\leq \sqrt{2}\RiskS[\wh{w}_{H_*}|\X_S]\leq\sqrt{2}\RiskS[\wh{w}_{H_*}|\bR^n],\] which combines with the above inequalities to imply that
\[
\Risk[\widehat{w}_{H_*}|\X_S]\leq \sqrt{2} \Risk[\widehat{w}_{H^*}|\X_S].
\]
However, the estimate $\widehat{w}_{H_*}$ cannot be seen as {\em adaptive} over the family of ``coaxial'' ellipsoids $\X_S^\kappa=\{x\in \bR^n:\,x^TSx\leq \kappa\}$, $\kappa\in K\subset \bR_+$, see, e.g., \cite{Lepski1991}. For instance, the maximal over $\X_S^\kappa$ risk
$\Risk[\widehat{w}_{H_*}|\X_S^\kappa]$ does not scale correctly for $\kappa\ll 1$ and $\kappa\gg 1$.

\subsubsection{Numerical illustration}\label{sectnumill}

In the above considerations, we treated matrix $S$ as part of the data. In fact, we can make $S$ a variable restricted to reside in a given computationally tractable convex subset $\S$ of the positive semidefinite cone, and look for minimal, over linear estimates {\sl and matrices $S\in\S$}, $S$-risk.
This can be done as follows. We consider a parametric family of problems with  $\tau$ in \rf{designproblem} being a parameter rather than a variable, and $S$ being a variable  restricted to reside in $\S$; then we apply bisection in $\tau$ to find the smallest value of $\tau$ for which the problem is feasible. With $S$ and linear estimate yielded by this procedure, the $S$-risk of the estimate clearly possesses near-optimality properties completely similar to those we have just established for the case of fixed $S$.
\par
As an illustration of these ideas, consider the following experiment. Let $[r;v]$ be {state of} pendulum with friction -- {the 2-dimensional } continuous time dynamical system {obeying the equations}
$$
\begin{array}{rcl}
\dot{{r}}&=&v,\\
\dot{v}&=&-{\nu}^2{r}-\kappa v+{w},\\
\end{array}
$$
where ${w}$ is the external input. Assuming this input constant on consecutive time intervals of duration $\Delta$, the sequence $z_\tau=[{r}(\tau\Delta);v(\tau\Delta)]$, $\tau=0,1,...$, obeys finite-difference equation
$$
z_{\tau}=Pz_{\tau-1}+Q{w}_\tau,\,\,\,\tau=1,2,...
$$
with
\[
P=\exp\Big\{\Delta\hbox{\scriptsize$\overbrace{\left[\begin{array}{cc}0&1\cr -{\nu}^2&-\kappa\cr\end{array}\right]}^{\vartheta}$}\Big\}, \;\;Q=\int_0^\Delta\exp\{s{\vartheta}\}\left[\begin{array}{c}0\\1\end{array}\right] ds;
\] here ${w}_\tau$ is the value of ${w}(\cdot)$ on the (continuous-time) interval $((\tau-1)\Delta,\tau\Delta)$. Assume that we are observing corrupted by noise positions ${r}_\tau={r}(\tau\Delta)$ of the pendulum on {the discrete-}time horizon $1\leq\tau\leq T$ and want to recover the inputs ${w}_s$, $T-K+1\leq s\leq T$. Denoting by $x=[z_0;{w}_1;{w}_2;...;{w}_T]$ the ``signal'' underlying our observations, we can easily build {a} $T\times(T+2)$ matrix $A$ and $1\times(T+2)$ matrices $B_t$ such that the trajectory ${r}:=[{r}_1;...;{r}_T]$ of pendulum's positions is given by ${r}=Ax$, and ${w}_t=B_tx$. {Given} noisy observations
$$
{\omega}=Ax+\sigma\xi,\;\xi\sim\N(0,I_T)
$$
of pendulum's (discrete time) trajectory, {we want to recover} inputs ${w}_t$, $1\leq t\leq T$, and their collections ${w}^K=[{w}_{T-K+1};{w}_{T-K+2};...;{w}_T]=B^{(K)}x$.\footnote{Note that estimating $w^K$ is not the same as ``standalone'' estimation of  {each individual} entry in ${w}^K$.}
\par
\begin{figure}[ht]
\centerline{
\begin{tabular}{cc}
\includegraphics[scale=0.35]{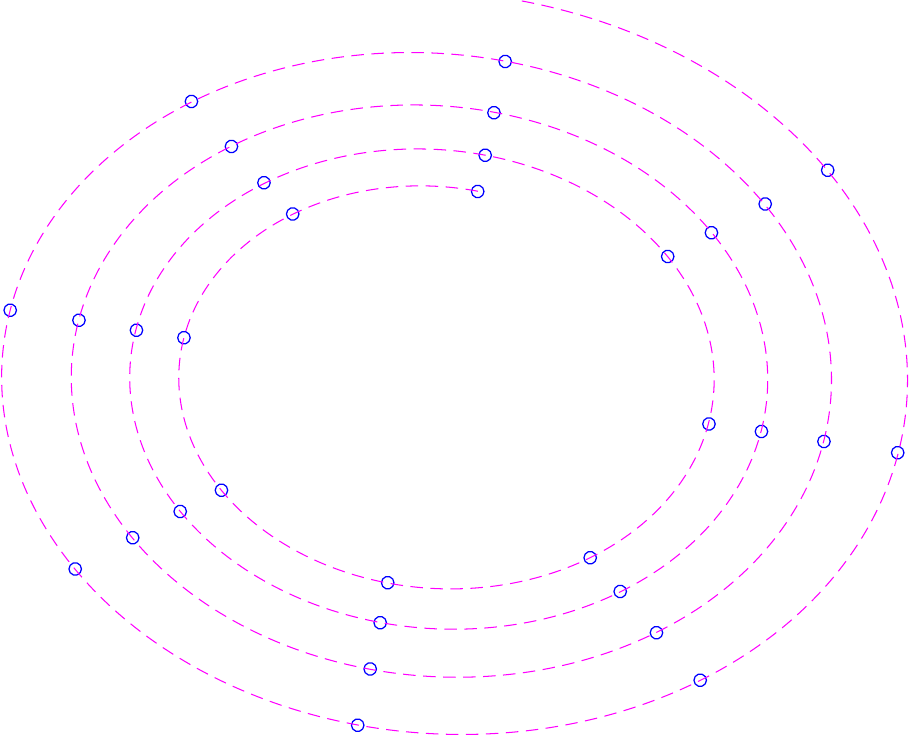}&\includegraphics[scale=0.35]{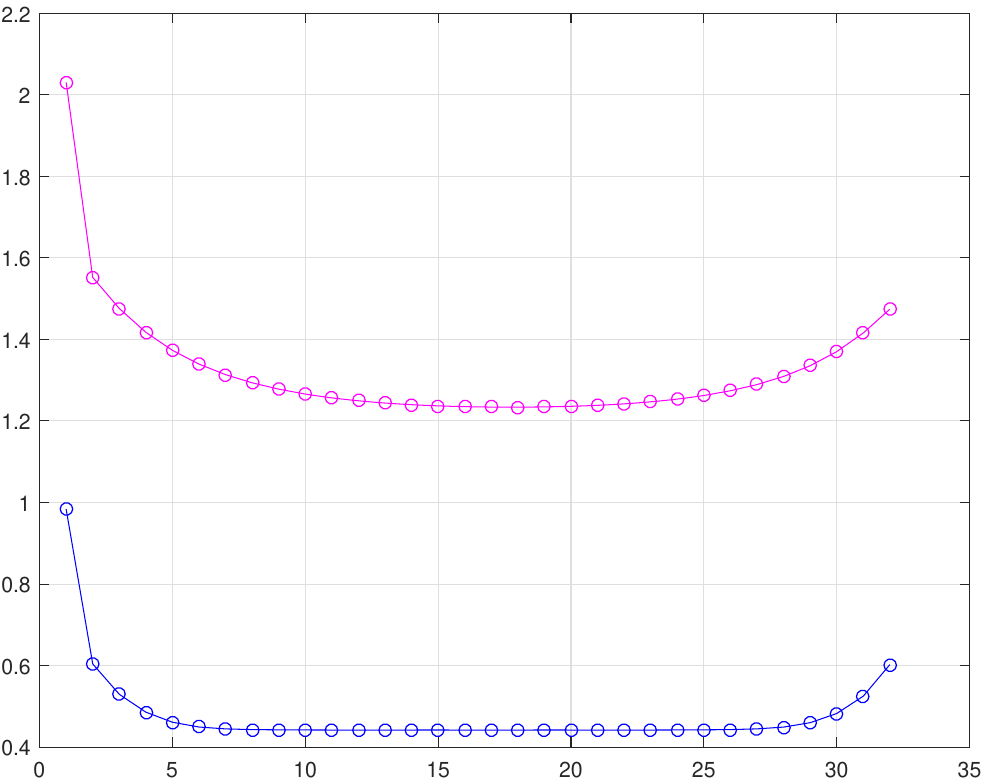}
\\(a)&(b)\\
\includegraphics[scale=0.35]{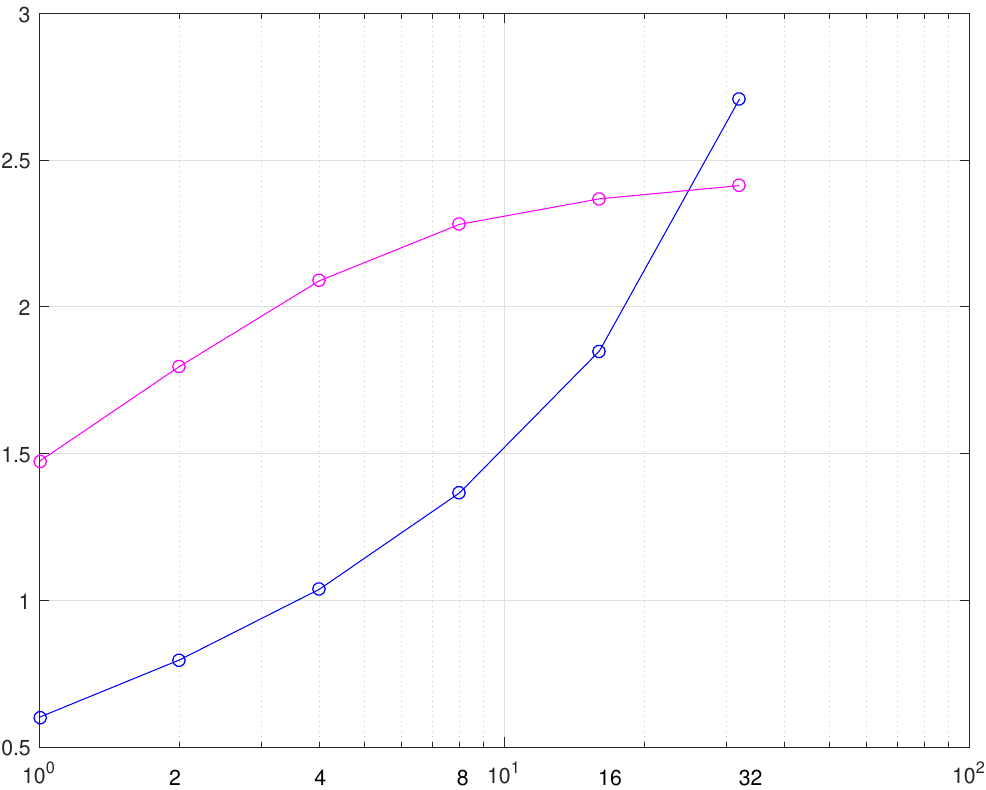}&\includegraphics[scale=0.35]{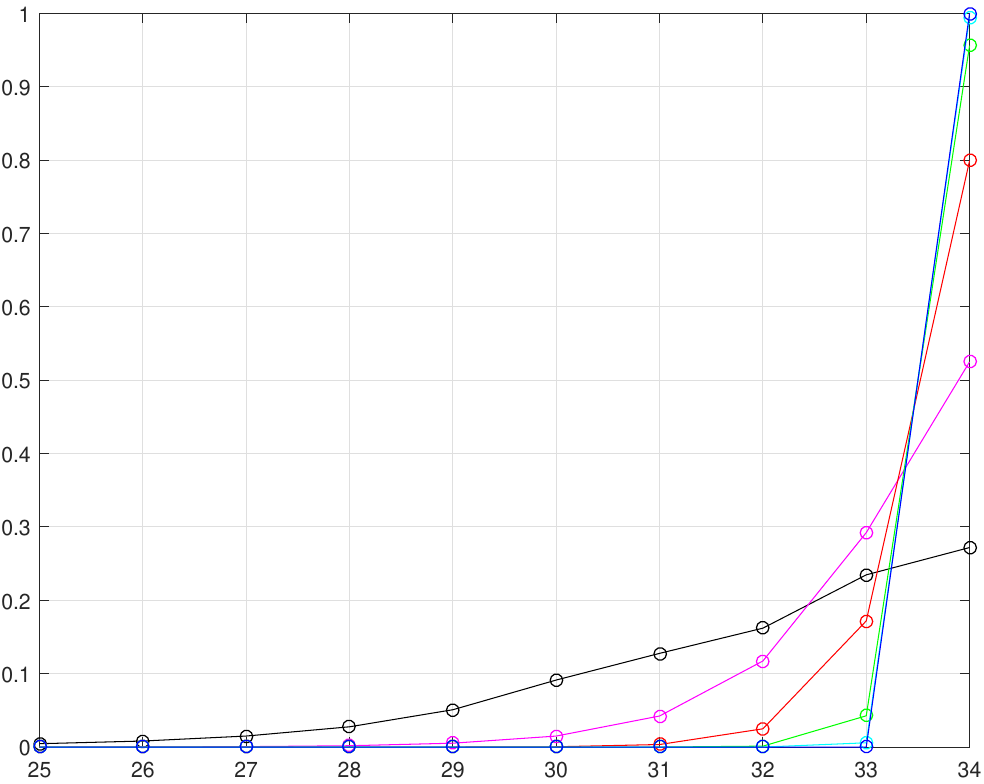}
\\(c)&(d)
\end{tabular}}
%
\caption{\label{pendulum} {\small Numerical illustration, section \ref{sectnumill}. (a): free motion (${w}\equiv0$) of pendulum in $({r},v)$-plane in continuous (dashed line) and discrete (circles) time. (b): Bayesian (blue) and worst-case (magenta) risks of recovering ${w}_t$ vs. $t=1,2,...,32$. (c): Bayesian (blue) and worst-case (magenta) risks of recovering ${w}^K:=[{w}_{T-K+1};{w}_{T-K+2};...;{w}_T]$ vs. $K$. (d): eigenvalues $\lambda_i(S_K)$ of $S_K$ ($K=32$ -- black, $K=16$ -- magenta, $K=8$ -- red, $K=4$ -- green, $K=2$ -- cyan, $K=1$ -- blue); we plot 10 largest eigenvalues of the $S$-matrices; the preceding 24 eigenvalues for all these matrices vanish.}}
\end{figure}

 We intend to process  our estimation problems by building the best, in terms of its $S$-risk taken over the entire space $\bR^{T+2}$ of signals, estimate; in our design, $S$ is not fixed in advance, {but is instead} restricted to be positive semidefinite with trace $\leq 1$. Thus, the problems we want to solve are of the form (cf. (\ref{designproblementirespace}))
\begin{equation}\label{examproblem}
\Opt[B]=\min_{\tau,H,S}\left\{\begin{array}{l}\tau:\left[\begin{array}{cc}\tau S&B^T-A^TH\cr
B-H^TA&I_T\cr\end{array}\right]\succeq0,\\~~~~~\sigma^2\Tr(H^TH)\leq\tau,\,S\succeq0,\,\Tr(S)\leq 1\end{array}\right\},
\end{equation}
where $B$ depends on what we want to recover ($B=B_t$ when recovering ${w}_t$, and $B=B^{(K)}$ when recovering ${w}^K$).
By Proposition \ref{prop6666}, the linear estimate $ H_{B,*}^T{\omega}$ yielded by an optimal solution $(\Opt[B],H_{B,*},S_{B,*})$ to the above (clearly solvable) problem is minimax optimal in terms of its $S$-risk $\RiskS[\cdot|\bR^{T+2}]$ {\sl taken with respect to $S=S_{B,*}$}, and the corresponding minimax optimal risk is exactly $\sqrt{\Opt[B]}$.
\begin{quote}
The rationale behind restricting $S$ to have its trace $\leq1$ is as follows. Imagine that we have reasons to believe that the entries in $x$ ``are of order of 1;'' the simplest way to model this belief is to assume that $x$ is uniformly distributed over the sphere $\S$ of radius $\sqrt{\dim x}=\sqrt{T+2}$. Under this assumption, the claim that an estimate $\widehat{w}(\cdot)$ has $S$-risk, taken over the entire space w.r.t. a matrix $S\succeq0$ with $\Tr(S)\leq1$, at most $\sqrt{\tau}$ means that
$$
\bE_{\xi\sim\N(0,I_T)}\{\|\widehat{w}(Ax+\sigma\xi)-B^Kx\|_2^2\}\leq \tau(1+x^TSx)\,\,\forall x.
$$
This relation, after taking expectation over the uniformly distributed over $\S$ signal $x$, implies that the expectation, over both $\xi$ and $x$, of the
 squared recovery risk is at most $2\tau$. Thus, optimising {the} $S$-risk over the linear estimates {\sl and} $S\succeq0$, $\Tr(S)\leq1$, can be interpreted as
safe minimization of the Bayesian risk taken w.r.t. a specific Bayesian prior (uniform distribution on $\S$). In this context, ``safety'' {means} that along with guarantees on the Bayesian risk, we get some meaningful upper bound on the expected $\|\cdot\|_2^2$-error of recovery applicable to {\sl every individual} signal.
\par
In view of the above considerations, { with some terminology abuse, below we} refer to the optimal value of (\ref{examproblem}) as to the {\sl Bayesian risk} of recovering $Bx$.
\end{quote}
In the experiment we are about to report, we use $\Delta=1$, $\kappa=0.05$ and select ${\nu}$ to make the eigenfrequency of the pendulum equal to 1/8; free motion of the pendulum in the $({r},v)$-coordinates is shown on Figure \ref{pendulum}. We used $\sigma=0.075$, $T=32$, and solved problem (\ref{examproblem}) for several ``$B$-scenarios.'' The results are presented on Figure \ref{pendulum} (b) -- (d). Plots (b) and (c) show the bound $\sqrt{2\Opt[B]}$, see above, on the Bayesian risk  along with {the optimal value of the risk of optimal linear recovery of $Bx$ for signals $x$ from the ball $\X$ of radius $\sqrt{T+2}$, as given by the optimal values of} (\ref{eq555}) (blue). Plot (b) shows what happens when recovering individual inputs ($B=B_t$, $t=1,2,...,T$) and displays the risks as functions of $t$; plot (c) shows the risks of recovering blocks $u^K=B^{(K)}x$ of inputs as functions of $K=1,2,4,...,32$. Finally, plot (d) shows the eigenvalues of the $S$-components of optimal solutions to problems (\ref{examproblem}) with $B=B^{(K)}$.\footnote{With $B=B_t$, $S$-components of optimal solutions to (\ref{examproblem}) turn out to be of rank 1 for all $t$.}

\subsection{Adding robustness}
\label{sectrob}
In this section we address the situation where the data $A,B$ of problems \rf{eq555} and (\ref{designproblem}) is not known exactly, and we are looking for estimates which are robust w.r.t. the corresponding data uncertainties. We lose nothing when restricting ourselves with problem (\ref{designproblem}), since  \rf{eq555} is the particular case $S=0$ of (\ref{designproblem}), with {ellitope} ${\X}$ given by (\ref{setX}). We intend to focus on the simplest case of {\sl unstructured norm-bounded uncertainty}
\begin{equation}\label{equnc}
[A;B]:=\left[\begin{array}{c}A\cr B\cr\end{array}\right]\in{\cal U}_r=\big\{[A;B]=[A_*;B_*]+E^T\Delta F:\Delta\in\bR^{p\times q},\|\Delta\|\leq r\big\};
\end{equation}
here $A_*\in\bR^{m\times n}$, $B_*\in\bR^{\nu\times n}$ are given {\sl nominal data}, and $E\in\bR^{p\times (m+\nu)}$, $F\in\bR^{q\times n}$ are given matrices.{\footnote{Recall that in the case of $P\neq I$ we have to replace matrices $A$, $B$ and $S$ with $AP$, $BP$ and $P^TSP$, respectively, and modify the definition of $\U_r$ accordingly: namely, when $[A;B]$ runs through the set $\U_r$, $[AP;BP]$ runs through
\[
\overline{\cal U}_r=\big\{[A;B]=[A_*P;B_*P]+E^T\Delta FP:\Delta\in\bR^{p\times q},\|\Delta\|\leq r\big\};
\]
where $A_*$, $B_*$ $E$ and $F$ are as in\rf{equnc}.} }
Our goal is to solve the {\sl robust counterpart}
{\begin{equation}\label{rcdesignproblem}
\begin{array}{l}
\hbox{\rm RobOpt}=\min\limits_{\tau,H,\lambda}\bigg\{\tau:\;\left[\begin{array}{cc}\sum_k\lambda_k S_k+\tau S&B^T-A^TH\cr
B-H^TA&I_\nu\cr\end{array}\right]\succeq0,\;\forall[A;B]\in{\cal U}\\
\multicolumn{1}{r}{\sigma^2\Tr(H^TH)+\phi_{\T}(\lambda)\leq\tau,\lambda\geq0}
\bigg\}
\end{array}
\end{equation}}
of problem (\ref{designproblem}).
 Plugging into (\ref{rcdesignproblem}) the parametrization of {$[A;B]$} via $\Delta$, the uncertainty-affected semidefinite constraint becomes
{\begin{eqnarray}
&&M(\lambda,\tau,H)+{\cal E}^T[H]\Delta{\cal F}+{\cal F}^T\Delta^T{\cal E}[H]\succeq0\;\;\forall (\Delta:\|\Delta\|\leq r),\nn
&&M(\lambda,\tau,H)=\left[\begin{array}{cc}\sum_k\lambda_k S_k+\tau S&B_*^T-A_*^TH\cr
B_*-H^TA_*&I_\nu\cr\end{array}\right],\label{eqrcnew}\\
&&{\cal E}[H]=[0_{p\times n},E_B-E_AH], \;{\cal F}=[F,0_{q\times\nu}],\nonumber
\end{eqnarray}
}
where
$$
E=[E_A,E_B]
$$
is the partitioning of the $p\times(m+\nu)$-matrix $E$ into the blocks comprised by the first $m$ and the last $\nu$ columns.
A well-known result of \cite{Boydetal} (see also \cite[section 8.2.1]{RO}) states that when ${\cal F}\neq0$ (this is the only nontrivial case), the semi-infinite Linear Matrix Inequality in (\ref{eqrcnew})  holds true if and only if there exists $\mu$ such that
$$
\left[\begin{array}{cc}
M(\lambda,\tau,H)-r^2\mu {\cal F}^T{\cal F}&[{\cal E}[H]]^T\cr
{\cal E}[H]&\mu I_p\cr\end{array}\right]\succeq0.
$$
It follows that the semi-infinite convex problem  (\ref{rcdesignproblem}) is equivalent to the explicit convex program
\begin{equation}\label{RCExplicit}
\begin{array}{rcl}
\hbox{\rm RobOpt}&=&\min\limits_{\tau,H,\lambda,\mu}\bigg\{\tau:
\,{\cal G}(H,\mu,\tau)\succeq0,
\\
&&\multicolumn{1}{r}{~~~~~~~~~~~\sigma^2\Tr(H^TH)+\phi_{\T}(\lambda)\leq\tau,\lambda\geq0\bigg\}}
\end{array}
\end{equation}
where \[
{\cal G}(H,\mu,\tau)=
\left[\begin{array}{c|c|c}\sum_k\lambda_k S_k+\tau S-\mu r^2F^TF&B_*^T-A_*^TH&\cr
\hline
B_*-H^TA_*&I_\nu&E_B^T-H^TE_A^T\cr
\hline
&E_B-E_AH& \mu I_p\cr\end{array}\right].
\]
The $H$-component of optimal solution  to (\ref{RCExplicit}) yields robust w.r.t. uncertainty (\ref{equnc}) estimate $H^T\omega$ of $Bx$ via observation $Ax+\sigma\xi$, and the expected $\|\cdot\|_2^2$-error of this estimate does not exceed $\hbox{\rm RobOpt}$, whatever be $x\in {\X}$  and  $[A;B]\in{\cal U}$.
\subsection{Byproduct on semidefinite relaxation}\label{secbypro}
A byproduct of {our main} observation (section \ref{mainobs}) we are about to {present} has nothing to do with statistics; it {relates to} the quality of the standard semidefinite relaxation. Specifically, given a quadratic from $x^TCx$ and an {ellitope} ${\X}$ represented by (\ref{setX}), consider the problem
\begin{equation}\label{Optminor}
\Opt_*=\max_{x\in {\X}} x^TCx=\max_{y\in\brX} y^TP^TCPy.
\end{equation}
This problem can be NP-hard (this is already so when ${\X}$ is the unit box and $C$ is positive semidefinite); however, $\Opt$ admits an efficiently computable upper bound given by {\sl semidefinite relaxation} as follows: whenever $\lambda\geq0$ is such that
$$
P^TCP\preceq \sum_{k=1}^K\lambda_kS_k,
$$
for $y\in \brX$ we clearly have
$$
[Py]^TCPy\leq \sum_k\lambda_ky^TS_ky\leq\phi_\T(\lambda)
$$
due to the fact that the vector with the entries $y^TS_ky$, $1\leq k\leq K$, belongs to $\T$. As a result, the efficiently computable quantity
\begin{equation}\label{appeq2}
\Opt=\min_{\lambda}\left\{\phi_\T(\lambda): \lambda\geq0,P^TCP\preceq\sum_k\lambda_kS_k\right\}
\end{equation}
is an upper bound on $\Opt_*$. We have the following
\begin{proposition}\label{propbyproduct} Let $C$ be a symmetric $n\times n$ matrix and ${\X}$ be given by ellitopic representation {\rm (\ref{setX})}, and let $\Opt_*$ and $\Opt$ be given by
{\rm (\ref{Optminor}) and (\ref{appeq2})}. Then
\begin{equation}\label{appeq3}
{\Opt\over {2\ln K+2\sqrt{\ln K}+1}} \leq \Opt_{*} \leq \Opt.
\end{equation}
\end{proposition}


\numberwithin{equation}{section}
\appendix

\section{Proofs}\label{sect:Proofs}
\subsection{Proofs for Section \ref{sect1}}
\subsubsection{Proof of Lemmas \ref{lem:lower1} and \ref{lem:lower2}}\label{prooflemlowers}
Since Lemma \ref{lem:lower1} is the particular case $S=0$ of Lemma \ref{lem:lower2}, we prove here only the latter statement. Let
$\widehat{w}(\cdot)$ be an estimate of $w=Bx$, and let $R$ be its $S$-risk, so that
$$
\forall (x\in X) : \bE_{\xi\sim \N(0,I_m)}\{\|\widehat{w}(Ax+\sigma\xi)-Bx\|_2^2\}\leq R^2(1+x^TSx),
$$
see \rf{S-riskdef}. Our intention is to bound $R$ from below. Observe that $xx^T\in\Q$ when $x\in\X$, whence $\|Bx\|_2=\sqrt{\Tr(Bxx^TB^T)}\leq M_*$ for all $x\in \X$, see (\ref{blueM_*}). It follows that projecting the estimate onto the $\|\cdot\|_2$-ball of radius $M_*$ centered at the origin, we can only reduce the risk of the estimate, and for the projected estimate the risk is at most $M_*$. Consequently, we
 can assume w.l.o.g. that
\be
R\leq M_* \ \&\ \|\widehat{w}(\omega)\|_2\leq M_*\;\;\forall \omega\in \bR^m.
\ee{eq4}
When taking expectation with respect to the distribution of the Gaussian vector $[\eta,\xi]$ with independent $\xi\sim\N(0,I_m)$ and $\eta\sim\N(0,Q)$, and  taking into account \rf{eq4}, we have for any $\gamma>0$
\begin{equation}\label{eq:0001}
\begin{array}{rcl}
\varphi(Q)&\leq&\bE_{[\xi,\eta]}\left\{\|\widehat{w}(A\eta+\sigma\xi)-B\eta\|_2^2\right\}\hbox{[by (\ref{blueeq:gaussopt})]}\\
&=&\bE_{\eta}\left\{\bE_{\xi}\left\{\|\widehat{w}(A\eta+\sigma\xi)-B\eta\|_2^2\right\}\right\}\\
&=&
\bE_\eta\left\{\bE_{\xi}\left\{\|\widehat{w}(A\eta+\sigma\xi)-B\eta\|_2^2\right\}1_{\eta\in \X}\right\}\\
&&
+\bE_\eta\left\{\bE_{\xi}\left\{\|\widehat{w}(A\eta+\sigma\xi)-B\eta\|_2^2\right\}1_{\eta\notin\X}\right\}
\\
&\leq& R^2\bE_\eta\left\{(1+\eta^TS\eta)1_{\eta\in \X}\right\}
+\bE_\eta \left\{[M_*+\|B\eta\|_2]^21_{\eta\notin \X}\right\}\\
&\leq& R^2\bE_\eta\left\{(1+\eta^TS\eta)\right\}+\bE_\eta\left\{\left[(1+1/\gamma)M_*^2+(1+\gamma)\|B\eta\|_2^2\right]1_{\eta\notin \X}\right\}\\
&\leq& R^2(1+\Tr(QS))+\bigg[(1+1/\gamma)M_*^2\delta+(1+\gamma) \underbrace{\bE
\left\{\|B\eta\|_2^21_{\eta\notin \X}\right\}}_{\I}\bigg];\\
\end{array}
\end{equation}
recall that $\delta\leq 1/5$ is an upper bound on the probability for $\eta\sim\N(0,Q)$ {\sl not} to belong to $\X$.
Let us upper-bound $\I$. We can find an orthogonal $U$ such that the matrix $U^TQ^{1/2}B^TBQ^{1/2}U$ is diagonal and can represent $\eta\sim \N(0,Q)$ as $\eta=Q^{1/2}U\zeta$ with $\zeta\sim\N(0,I)$; denoting $\overline{Z}=\{\zeta:\,Q^{1/2}\zeta\not\in \X\}$, we get
\[
\I=\bE\left\{\zeta^T{U^T}Q^{1/2}B^TBQ^{1/2}{U}\zeta 1_{\zeta\in \overline{Z}}\right\},
\]
with
\[
\Prob\{\zeta\in \overline{Z}\}\leq \delta.
\]
Recalling that the matrix $U^TQ^{1/2}B^TBQ^{1/2}U$ is diagonal and $\succeq0$, we have
\[\begin{array}{l}
\I\leq\overbrace{\Tr(U^TQ^{1/2}B^TBQ^{1/2}U)}^{=\Tr(Q^{1/2}B^TBQ^{1/2})}\max_{1\leq i\leq n}\bE\{ \zeta_i^2 1_{\zeta\in \overline{Z}}\}
=
\|BQ^{1/2}\|_{2}^2\max_{1\leq i\leq n}{1\over \sqrt{2\pi}}\int_{-\infty}^\infty s^2\e^{-s^2/2}\chi_{\overline{Z}}(s)ds\\
\end{array}
\]
where $\chi_{\overline{Z}}(s)$ is the conditional, given that $\zeta_i=s$, probability for $\zeta\sim\N(0,I_{n})$ to belong to $\overline{Z}$, so that
$0\leq\chi_{\overline{Z}}(s)\leq 1$, and \[
{1\over \sqrt{2\pi}}\int_{-\infty}^\infty \chi_{\overline{Z}}(s) \e^{-s^2/2}ds\leq\delta.
 \] {We conclude (see  Lemma \ref{lem: referee} below) that}
$$
\I\leq \|BQ^{1/2}\|_2^2 \sqrt{{2\over\pi}}\int_{q_{1-{\delta/ 2}}}^\infty s^2\e^{-s^2/2}ds
$$
where $q_t$ is the $t$-quantile of the standard normal distribution:
$$
{1\over\sqrt{2\pi}}\int_{-\infty}^{q_t}\e^{-s^2/2}ds=t,\;0< t<1.
$$
On the other hand, for $\delta\leq 1/5$ one has
\[
\sqrt{{2\over\pi}}\int_{q_{1-{\delta/ 2}}}^\infty s^2\e^{-s^2/2}ds\leq 2 \sqrt{{2\over\pi}} q^2_{1-\delta/2}\int_{q_{1-{\delta/ 2}}}^\infty\e^{-s^2/2}ds=2 q^2_{1-\delta/2} \delta,
\]
and
\[
\I\leq 2\|BQ^{1/2}\|_2^2 q^2_{1-\delta/2}\delta.
\]
When substituting the latter bound into \rf{eq:0001}
we conclude that
\[
\varphi(Q)\leq R^2(1+\Tr(QS))+\bigg[(1+1/\gamma)M_*^2+2(1+\gamma)\|BQ^{1/2}\|_2^2 q^2_{1-\delta/2}\bigg]\delta
\]
Hence, when optimizing in $\gamma>0$
we obtain
\[
\varphi(Q)\leq R^2(1+\Tr(QS))+[M_*+\sqrt{2}q_{1-\delta/2}\|BQ^{1/2}\|_2]^2\delta.
\]
When passing to the limit as $R\to \RiskoptS[\X]+0$,
we come to
\[
\varphi(Q)\leq \Riskopt^2[\X](1+\Tr(QS))+[M_*+\sqrt{2} q_{1-\delta/2} \|BQ^{1/2}\|_{2}]^2\delta,
\]
what is \rf{eq:0002} for $S=0$. Finally, when $Q\in\rho\Q$, by \rf{blueM_*} we get
$\|BQ^{1/2}\|_{2}\leq\sqrt{\rho}M_*$,
and we arrive at \rf{eq:000222}.
\hfill\qed
{\begin{lemma}\label{lem: referee}
Let $f:\,\bR\to [0,1]$ such that ${1\over \sqrt{2\pi}}\int_{-\infty}^\infty f(s)\exp\{-s^2/2\}ds\leq\delta$. Then  \[
\int s^2f(s)\exp\{-s^2/2\}ds\leq
\int_{|s|\geq q_{1-\delta/2}}s^2\exp\{-s^2/2\}ds.
\]
\end{lemma}
{\bf Proof.} Indeed,  let us denote
\[
\bar{f}(s)=\left\{\begin{array}{ll}0,&|s|<q_{1-\delta/2}\\
1,&|s|\geq q_{1-\delta/2}.\\
\end{array}\right.
\]
We have
\bse
\bar{f}(s)-f(s)&\leq&0,\;\mbox{for}\;|s|\leq q_{1-\delta/2},\\
\bar{f}(s)-f(s)&\geq&0,\;\mbox{for}\;|s|\geq q_{1-\delta/2},
\ese and
\[
 \int_{-\infty}^\infty f(s)\exp\{-s^2/2\}ds\leq\sqrt{2\pi}\delta=\int_{-\infty}^\infty
\bar{f}(s)\exp\{-s^2/2\}ds.
\]
Thus,
\bse
\lefteqn{\int_{-\infty}^\infty
s^2\bar{f}(s)\exp\{-s^2/2\}ds-\int_{-\infty}^\infty
s^2f(s)\exp\{-s^2/2\}ds}\\
&=&\int_{|s|\geq q_{1-\delta/2}}s^2\underbrace{[\bar{f}(s)-f(s)]}_{\geq0\hbox{\ \tiny when $|s|\geq q_{1-\delta/2}$}}\exp\{-s^2/2\}ds-\int_{|s|<q_{1-\delta/2}}s^2\underbrace{[f(s)-\bar{f}(s)]}_{\geq0\hbox{\
\tiny when $|s|<q_{1-\delta/2}$}}\exp\{-s^2/2\}ds\\
&\geq& q_{1-\delta/2}^2\int_{|s|\geq q_{1-\delta/2}}[\bar{f}(s)-f(s)]\exp\{-s^2/2\}ds - q_{1-\delta/2}^2\int_{|s|<q_{1-\delta/2}} [f(s)-\bar{f}(s)]\exp\{-s^2/2\}ds\\
&=&q_{1-\delta/2}^2\left[\int_{-\infty}^\infty
[\bar{f}(s)-f(s)]\exp\{-s^2\}ds\right]\geq0.~~~~~~~~~~~~~~~~~~~~~~~~~~~~~~~~~~~~~~~~~~~~~~~~~~~~~~~~~\mbox{\qed}
\ese
}

\subsubsection{Proof of Lemma \ref{Opto}}
We set (cf. \rf{conebT})
\[
\bT=\cl \{[t;\tau]\in\bR^K\times\bR:\tau>0,\tau^{-1}t\in\T\}\subset \bR^{K+1}_+;
\]
recall that $\bT$ is a closed and pointed convex cone in $\bR^{K+1}$ with a nonempty interior such that
\[
\T=\{t:\exists \tau\leq 1: [t;\tau]\in\bT\}=\{t:[t;1]\in\bT\}.
\]
Note that \rf{starrho} is nothing but the conic problem
\begin{equation}\label{eq:nowdiez}
\Opt_*=\max\limits_{Q,G,t}\left\{\Tr(\brB Q\brB^T)-\Tr(G): \begin{array}{l}\left[\begin{array}{cc}G&\brB Q\brA^T\cr \brA Q\brB^T&\sigma^2I_m+\brA Q\brA^T\cr\end{array}\right]\succeq0\\
Q\succeq0,\,[t;1]\in\bT,\,
\Tr(QS_k)\leq t_k,\,1\leq k\leq K
\end{array}\right\}.
\end{equation}
This problem clearly is strictly feasible (since $\inter  \T$ contains a positive vector) and bounded (the latter is due to $\sum_kS_k\succ0$), so that its optimal value is equal to the optimal value of its conic dual problem, and all we need in order to prove (\ref{miracle}) is to verify is that the latter problem is equivalent to (\ref{eq555}).\par
Let us build the dual to \rf{eq:nowdiez} (for ``guidelines,'' see Appendix \ref{conicd}). Note that
the cone dual to $\bT$ is
\[
\bT_*=\{[g;s]: s\geq \phi_{\T}(-g)\}.
\]
Denoting the Lagrange multiplier for the first $\succeq$-constraint in \rf{eq:nowdiez} by $\left[\begin{array}{cc}U&V\cr V^T&W\cr\end{array}\right]\succeq0$, for the second $\succeq$-constraint by $L\succeq0$, for $\leq$-constraints by $-\lambda$,  $\lambda\in\bR^K_+$, and for the constraint $[t;1]\in\bT$ -- by $[g;s]\in\bT_*$, multiplying the constraints by the multipliers and summing up the results, we see that
the constraints in \rf{eq:nowdiez} imply that on the feasible set of \rf{eq:nowdiez} it holds
\begin{equation}\label{eq444}
\begin{array}{l}
-\Tr(UG)-\Tr(Q[\brB^TV\brA+\brA^TV^T\brB])-\Tr(Q[\brA^TW\brA])-\Tr(LQ)+\sum\limits_k\lambda_k\Tr(QS_k)\\-\sum\limits_k\lambda_kt_k
-\sum\limits_kg_k t_k\leq \sigma^2\Tr(W)+s.
\end{array}
\end{equation}
Now to get the dual to \rf{eq:nowdiez} problem, we need to impose on the Lagrange multipliers the constraint that the left hand side in (\ref{eq444}) is identically in $Q,G,t$ equal to the objective $\Tr(\brB Q\brB^T)-\Tr(G)$ of \rf{eq:nowdiez}, and to minimize over the multipliers under this constraint (in addition to those introduced when specifying the multipliers) the right hand side of (\ref{eq444}). Thus, the problem dual to \rf{eq:nowdiez} is
\bse
[\Opt_*&=&]\min\limits_{U,V,W,L,\lambda,g,s}\bigg\{\sigma^2\Tr(W)+s:\;
\left[\begin{array}{cc}U&V\cr V^T&W\cr\end{array}\right]\succeq0,\,
L\succeq0,\,\lambda\geq0,\,s\geq \phi_{\T}(-g),\\
&&g_k=-\lambda_k,\,1\leq k\leq K,\,U=I_\nu,\,-\brB^TV\brA-\brA^TV^T\brB-\brA^TW\brA-L+\sum_k\lambda_k S_k=\brB^T\brB
\bigg\}\\
&=&\min\limits_{V,W,\lambda,s}\left\{\sigma^2\Tr(W)+s: \;\begin{array}{l}W\succeq V^TV,\,\lambda\geq0,\,s\geq \phi_\T(\lambda),\\
\sum_k\lambda_k S_k\succeq \brB^T\brB+\brB^TV\brA+\brA^TV^T\brB+\brA^TW\brA\\
\end{array}\right\}\\
&=&\min\limits_{V,W,\lambda}\left\{\sigma^2\Tr(V^TV)+\phi_\T(\lambda):\;\begin{array}{l} W=V^TV,\,\lambda\geq0\\
\sum_k\lambda_k S_k\succeq \brB^T\brB+\brB^TV\brA+\brA^TV^T\brB+\brA^TW\brA\\
\end{array}\right\}\\
&=&\min\limits_{V,\lambda}\left\{\sigma^2\Tr(V^TV)+\phi_\T(\lambda):\; \sum_i\lambda_k S_k\succeq (\brB+V\brA)^T(\brB+V\brA),\,\lambda\geq0\right\},
\ese
that is,
$\Opt_*=\Opt$ (substitute $H=-V^T$ in (\ref{eq555})).\hfill\qed
\subsubsection{Proof of Lemma \ref{lem1}}
Representing $\eta=Q^{1/2}\zeta$ with $\zeta\sim\N(0,I_n)$, we reduce the situation to the one where $(Q,S)$ is replaced with $(I_{n},\bar{S}=Q^{1/2}SQ^{1/2})$, so that it suffices to prove (\ref{then1}) in the special case of $Q=I_{n}$. Moreover, we clearly can assume that $S$ is diagonal with diagonal entries $s_i\geq0$, $1\leq i\leq n$, so that $\rho=\sum_is_i$. Now the relation we should prove reads
$$
\Prob_{\eta\sim\N(0,I_{n})} \left\{\sum_{i=1}^{n}s_i\eta_i^2>1\right\}\leq \e^{-{1-\rho+\rho\ln(\rho)\over 2\rho}}.
$$
Let $\gamma\geq0$ be such that $2\gamma\max_is_i<1$. Then
\[
\begin{array}{l}
\ln\left(\bE_\eta\{\exp\{\gamma \sum_{i=1}^{n}s_i\eta_i^2\} \}\right)=\sum_{i=1}^{n}\ln\left(\bE_\eta\{\exp\{\gamma s_i\eta_i^2\}\}\right)=-{1\over 2}\sum_{i=1}^{n}\ln(1-2\gamma s_i),
\end{array}
\]
what implies the first inequality of \rf{then1}. Furthermore,
 for $0\leq\gamma<{1\over 2\max_i s_i}\leq {1\over 2\rho}$,
$$
\ln\left(\bE_\eta\left\{\exp\left[\gamma \sum_{i=1}^{n}s_i\eta_i^2\right] \right\}\right)\leq -{1\over 2}\ln(1-2\gamma\rho)
$$
(indeed, the convex function $-{1\over 2}\sum_{i=1}^{n}\ln(1-2\gamma s_i)$ of $s$ varying in the simplex $\{s\geq0,\sum_is_i=\rho\}$ attains its maximum at a vertex of the simplex). Specifying $\gamma={1-\rho\over 2\rho}$, we conclude that
\bse
\Prob\left\{\sum_{i=1}^{n}s_i\eta_i^2>1\right\}&\leq& \bE_\eta\{\exp\{\gamma \sum_{i=1}^{n}s_i\eta_i^2\} \}\exp\{-\gamma\}
\leq \exp\{-{1\over 2}\ln(1-2\gamma\rho)-\gamma\}\\
&=&\exp\{-{1-\rho+\rho\ln(\rho)\over 2\rho}\},
\ese
as claimed. \qed
\subsubsection{Proof of Corollary \ref{simplefactor}}
Observe that ${\X}$ contains a point $\bar{x}$ with $$\|\bar{x}\|_2\geq r:={\sqrt{T}\over \cond(\T)\sqrt{\varkappa}}.$$
Indeed, by definition of $\cond(\T)$, $\T$
contains a vector $\bar{t}$ with all entries $\geq T/\cond^2(\T)$; let now $\bar{x}=re$, where $e$ is the eigenvector of the matrix $S=\sum_{k=1}^KS_k$ corresponding to the minimal eigenvalue $\varkappa$ of this matrix. We have (recall that $S_k\succeq 0,\, k=1,...,K$)
\[
\bar{x}^TS_k\bar{x}\leq \varkappa r^2
=T/{\cond}^2(\T)\leq \bar{t}_k,\;\;1\leq k\leq K,
\]
that is, $\bar{x}\in {\X}$. Selecting  the largest $t\in[0,1]$ such that
$t\|A\bar{x}\|_2\leq c\sigma$, where $c$ is the positive absolute constant from \rf{eq:slower1}, we conclude  by \rf{eq:slower1} that $\Riskopt[\X]\geq t\|B\bar{x}\|_2$, or
\[
\Riskopt[\X]\geq \|B\bar{x}\|_2\min\left[1,{c\sigma\over \|A\bar{x}\|_2}\right].
\]
Hence, we get
\bse\Riskopt[{\X}]&\geq& \sigma_{\min(B)}\|\bar{x}\|_2\min\left[1, {c\sigma\over \|A\|\,\|\bar{x}\|_2}\right]
\geq {\|B\|\over \cond(B)} \min\left[\|\bar{x}\|_2,{\sigma\over \|A\|}\right]\\
&\geq& {\|B\|\over \cond(B)} \min\left[r,{\sigma\over \|A\|}\right]=
{\|B\|\over \cond(B)} \min\left[{\sqrt{T}\over \cond(\T)\sqrt{\varkappa}},{\sigma\over \|A\|}\right]
\ese
Note that the quantity $M_*=\max_{Q\in \Q}\|BQ^{1/2}\|_2$ admits simple bound:
\[
M_*\leq \|B\|\sqrt{T/\varkappa}
\]
(indeed, since $\sum_{k}\Tr(QS_k)\leq T$ for all $Q\in \Q$, one has $\Tr(Q\sum_kS_k)\leq T$, whence $\Tr(Q)\leq T/\varkappa$
by the origin of $\varkappa$, and therefore
$M_*^2=\Tr(BQB^T)\leq \|B^TB\|\Tr(Q)\leq\|B\|^2T/\varkappa$).
As a result,
\bse
{M_*\sqrt{K}\over \tRiskopt[{\X}]}&\leq&
c'\cond(B)\sqrt{TK\over \varkappa}\max\left[\cond(\T)\sqrt{\varkappa\over T},\,{ \|A\|\over\sigma}
\right]\leq
c'\cond(B)\sqrt{K}\left[\cond(\T)+{\|A\| \sqrt{T}\over\sigma\sqrt{\varkappa}}\right]\\
\ese
with an absolute constant $c'$; together with \rf{nearoptimal} this implies \rf{eq:simplef}.\hfill\qed
\subsection{Proofs for Section \ref{sectextensions}}
\subsubsection{Proof of Lemma \ref{OptSo}}
\paragraph{1$^o$}We claim that (\ref{eq:diez}) is a strictly feasible conic problem with bounded level sets of the objective (the sets where the objective is $\geq a$, for every fixed $a\in\bR$); in particular, the problem is solvable.
\par
Indeed, strict feasibility follows from the fact that the interior of the cone $\bT$ contains a positive vector, see assumptions on $\T$ in section  \ref{sitgoal}.  Further, the projections of the feasible set onto the $[v;s]$- and $W$-spaces are bounded (the first -- since at a feasible solution it holds $0\leq s\leq1$, and the second -- due to the boundedness of the set of $v$-components of feasible solutions combined with $\sum_kS_k\succ0$). Boundedness of a level set of the objective follows from the fact that if a sequence of feasible solutions
$\{(W_i,G_i,[v^i;s^i]),i=1,2,...\}$ goes to $\infty$, then, by the above, the sequence $\{W_i,[v^i;s^i]\}$ is bounded, so that
 $\|G_i\|\to\infty$ as $i\to\infty$; since $G_i\succeq0$ due to the constraints of the problem, we have $\Tr(G_i)\to\infty$ as $i\to\infty$, which combines with boundedness of $\{W_i\}$ to imply that the objective along our sequence of feasible solutions goes to $-\infty$, which is impossible for a sequence of feasible solutions from a level set of the objective.
\paragraph{2$^o$} Our next claim is that at an optimal solution $(W,G,[v;s])$ to \rf{eq:diez} one has $s>0$.
\par Indeed, otherwise $v=0$ due to $[v;s]\in\bT$ and the origin of $\bT$, whence $W=0$ due to $W\succeq0$ and $\sum_kS_k\succ0$; besides this, $G\succeq0$, so that assuming $s=0$, we see that $\Opt_*=0$, which clearly is not the case: $\bT$ contains a vector $[\bar{v};\bar{s}]$ with, say, $\bar{s}=0.1$ and positive $\bar{v}$, implying that for some $\bar{\tau}>0$ and all $\tau\in[0,\bar{\tau}]$ tuples
$$
W_\tau=\tau I,G_\tau=[\sigma^2\bar{s}]^{-1}[BW_\tau A^TAW_\tau B^T]=[\sigma^2\bar{s}]^{-1}\tau^2BA^TAB^T,[\bar{v};\bar{s}]
$$
are feasible solutions to (\ref{eq:diez}); since $B\neq0$, for small positive $\tau$ the value of the objective of (\ref{eq:diez}) at such a solution is positive, which would be impossible when $\Opt_*=0$.
\par
{Furthermore,}
 observe that if  $(W,G,v,s)$ is an optimal solution to \rf{eq:diez} (whence, as we already know, $s>0$), when replacing $G$ with the matrix
$$
\bar{G}:=BWA^T(\sigma^2s I_m +A WA^T)^{-1}A WB^T
$$
(so that $G\succeq \bar{G}$ and $(W,\bar{G},t,s)$ is feasible for \rf{eq:diez}),
we keep the solution optimal, thus
\[
\Opt_*=\Tr\left(B [W-WA^T(\sigma^2s I_m +A WA^T)^{-1}A W]B^T\right).
\]
\paragraph{3$^o$} {To complete the proof of the lemma it suffices to show}
that the conic dual to \rf{eq:diez} is equivalent to (\ref{designproblem}); since \rf{eq:diez}, as we have already mentioned, is strictly feasible and bounded, this would imply that
$\Opt=\Opt_*.$
\par
To build the problem dual to \rf{eq:diez}, let the Lagrange
multipliers for the constraints be, respectively, $\left[\begin{array}{cc}U&V\cr V^T&Z\cr\end{array}\right]\succeq0$, $L\succeq0$, $-\lambda$, $\lambda\in\bR^K_+$, $-\tau$, $\tau\geq0$, and $[g;r]\in\bT_*$,
where
$$
\bT_*=\{[g;r]:\; r\geq \phi_{\T}(-g)\}
$$
is the cone dual to $\bT$. Taking inner products of the constraints of \rf{eq:diez} with the multipliers and summing up the results, we arrive at the aggregated constraint
$$
\begin{array}{l}
\Tr(GU) +\Tr(W[A^TV^TB+B^TVA+A^TZA+L-\sum_k\lambda_k S_k-\tau S])\\
\multicolumn{1}{r}{+\sum_k[\lambda_k+g_k]v_k+s[\sigma^2\Tr(Z)-\tau+r]+\tau\geq0}\\
\end{array}
$$
To get the dual problem, we impose on the multipliers the restriction for the resulting inequality to have the homogeneous in $W,G,v,s$ component identically equal to {\sl minus} the objective of \rf{eq:diez}, which amounts  to the relations
\bse
&&U=I_\nu,\; \tau=r+\sigma^2\Tr(Z),\;g_k=-\lambda_k\;\forall k,\\
&&[A^TV^TB+B^TVA+A^TZA+L-\sum_k\lambda_k S_k-\tau S]=-B^TB.
\ese
Under these relations, the aggregated constraint reads
$$
\Tr(BWB^T-G)\leq \tau
$$
for all feasible solutions to \rf{eq:diez}, thus $\Opt_*\leq\tau$. Therefore, the problem dual to \rf{eq:diez} is to minimize the resulting upper bound on $\Opt_*$, that is, the dual is
$$
\min_{\tau,V,Z,L,\lambda,[g;r]}\left\{\tau:\begin{array}{l}\left[\begin{array}{cc}
I_\nu&V\cr V^T&Z\cr\end{array}\right]\succeq0,\,
L\succeq0,\,\lambda\geq0,\,\tau\geq0,\,r\geq\phi_{\T}(-g)\,\\
B^TB+A^TV^T+VA+A^TZA=\sum_k\lambda_k S_k+\tau S-L\\
g=-\lambda,\tau=r+\sigma^2\Tr(Z)\\
\end{array}\right\}.
$$
Now partial minimization in $Z$ and $r$ results in $Z=V^TV$ which, after eliminating $L$ and $[g;r]$, reduces the dual problem to
$$
\min_{\tau,V,\lambda}\left\{\tau:\begin{array}{l}(B+VA)^T(B+VA)\preceq\sum_k\lambda_k S_k+\tau S,\\
\lambda\geq0,\,\tau\geq \phi_{\T}(\lambda)+\sigma^2 \Tr(V^TV)\\
\end{array}\right\}.
$$
The resulting problem clearly is equivalent to (\ref{designproblem}) (substitute $V=-H^T$). Thus, (\ref{equal}) is proved.\hfill\qed
\subsubsection{Proof of Proposition \ref{prop6666}}
Under the premise of the proposition, the feasible set of (\ref{designproblementirespace}) is nonempty, and the objective clearly goes to $\infty$ along every going to $\infty$ sequence of feasible solutions $(\tau_i,H_i)$, implying that the problem is solvable. The optimal value $\Opt$ in the problem clearly is positive due to $\sigma>0$ and $B\neq0$. Now assume that (\ref{minmaxoptimality}) does not hold, so that there exists $\alpha$ and estimate $\widehat{w}_*(\cdot)$ such that
\begin{equation}\label{contradiction}
\alpha<\Opt\; \&\; \bE_{\xi\sim\N(0,I_m)}\{\|\widehat{w}_*(Ax+\sigma\xi)-Bx\|_2^2\}\leq \alpha(1+x^TSx)\,\,\forall x\in\bR^n,
\end{equation}
and let us lead this assumption to contradiction.
\par
Consider the conic problem (cf. (\ref{eq:diez}))
\begin{equation}\label{eq:diezunb}
\Opt_*=\max_{W,G,s}\left\{\Tr(BWB^T)-\Tr(G):\begin{array}{l}
\left[\begin{array}{cc}G&BWA^T\cr A WB^T&\sigma^2sI_m+A WA^T\cr\end{array}\right]\succeq0,\\
W\succeq0,\,\Tr(WS)+s\leq1,\,s\geq0\\
\end{array}\right\}.
\end{equation}
This conic problem clearly is strictly feasible; the same argument as in the case of (\ref{eq:diez}) shows that the conic dual of this problem is equivalent to (\ref{designproblementirespace}) and therefore is feasible. By Conic Duality Theorem, it follows that both (\ref{eq:diezunb}) and (\ref{designproblementirespace}) have equal optimal values, and since $\sigma>0$, $B\neq0$, $\Opt$ is positive. Thus,
\[
\Opt_*=\Opt>0.
\]
This relation, due to $\alpha<\Opt$, implies that there is a feasible solution to (\ref{eq:diezunb}) with the value of the objective $>\alpha$. Since the problem is strictly feasible, feasible solutions with $s>0$ are dense in the feasible set, implying that the above feasible solution, let it be $(\widehat{W},G,\widehat{s})$, can be selected to have $\widehat{s}>0$. Further, keeping $\widehat{W}$ and $\widehat{s}$ intact and replacing $G$ with $\widehat{G}=B\widehat{W}A^T[\sigma^2\widehat{s}I_m+A \widehat{W}A^T]^{-1}A\widehat{W}B^T$, we preserve feasibility and can only increase the objective of (\ref{eq:diezunb}). The bottom line is that we can point out a feasible solution $(\widehat{W},\widehat{G},\widehat{s})$ to (\ref{eq:diezunb}) such that
\begin{equation}\label{goodsolution}
\begin{array}{l}
\widehat{\alpha}:=\Tr(B^T[\widehat{W}-\widehat{W}A^T[\sigma^2\widehat{s}I_m+A \widehat{W}A^T]^{-1}A\widehat{W}]B) >\alpha,\\
\widehat{s}>0,\; \widehat{W}\succeq 0,\;\Tr(\widehat{W}S)+\widehat{s}\leq 1.
\end{array}
\end{equation}
Observe that
\begin{equation}\label{varphiagain}
\widehat{\alpha}=\widehat{s}\varphi(\widehat{s}^{-1}\widehat{W})
\end{equation}
(see (\ref{diezrho})). Now let $\eta\sim \N(0,\widehat{s}^{-1}\widehat{W})$ be independent of $\xi\sim\N(0,I_m)$. We have
$$
\begin{array}{rcl}
\bE_{[\eta;\xi]}\{\|\widehat{w}_*(A\eta+\sigma\xi)-B\eta\|_2^2\}&=&\bE_{\eta}\left\{\bE_{\xi}\{\|\widehat{w}_*(A\eta+\sigma\xi)-B\eta\|_2^2\}\right\}\\
&\leq& \bE_{\eta}\left\{\alpha(1+\eta^TS\eta)\right\}\hbox{\ [by (\ref{contradiction})]}\\
&=&\alpha(1+\widehat{s}^{-1}\Tr(\widehat{W}S)).
\end{array}
$$
By (\ref{blueeq:gaussopt}), the initial quantity in this chain is $\geq\varphi(\widehat{s}^{-1}\widehat{W})=\widehat{s}^{-1}\widehat{\alpha}$ (see (\ref{varphiagain})), so that the chain yields $\widehat{s}^{-1}\widehat{\alpha}\leq \alpha(1+\widehat{s}^{-1}\Tr(\widehat{W}S))$, that is,
$$
\widehat{\alpha}\leq \alpha(\widehat{s}+\Tr(\widehat{W}S))\leq\alpha,
$$
where the last $\leq$ stems from the last inequality in (\ref{goodsolution}). The resulting inequality contradicts the first inequality in (\ref{goodsolution}); we have arrived at the desired contradiction.
\hfill\qed
\subsubsection{Proof of Proposition \ref{propbyproduct}}
We need the following
\begin{lemma}\label{largedev} Let $S$ be a positive semidefinite $\brn\times \brn$ matrix with unit trace and let $\xi$ be a Rademacher $\brn$-dimensional random vector (i.e., the entries in $\xi$ are independent and take values $\pm1$ with probabilities $1/2$). Then
{ for all $x\geq 0$ one has
\be
\Prob\left\{\xi^TS\xi\geq 1+2x+2\sqrt{x}\right\}\leq e^{-x}.
\ee{appeq1001}
}
\end{lemma}
{\bf Proof.} Let $S=\sum_{i=1}^{\brn}\lambda_ig_ig_i^T$ be the eigenvalue decomposition of $S$, so that $\lambda_i\geq0$, $\sum_i\lambda_i=1$ and $\|g_i\|_2=1$. Then {for $\kappa\geq 0$,}
$$
\bE\Big\{\exp\big\{{\kappa}\xi^TS\xi\big\}\Big\}=\bE\left\{\exp\Big\{{\kappa}\sum_i\lambda_i(g_i^T\xi)^2\Big\}\right\}
$$
is a convex function of $\lambda$ and therefore it attains its maximum over nonnegative vectors $\lambda$ with unit sum of entries at a basic orth. {
On the other hand, when $S=gg^T$ with unit vector $g$ one has for $0\leq \kappa<\half$ and $\eta\sim \N(0,1)$  independent of $\xi$:
$$
\begin{array}{l}
\bE\Big\{\exp\left\{\kappa(g^T\xi)^2\right\}\Big\}=\bE_\xi\Big\{\bE_{\eta}\Big\{\e^{\eta\sqrt{2\kappa} g^T\xi}\Big\}\Big\}=\bE_\eta\Big\{\bE_{\xi}\Big\{\e^{ \eta\sqrt{2\kappa} g^T\xi}\Big\}\Big\}=\bE_\eta\Big\{\prod_{i=1}^{\bar{n}}\bE_{\xi_i}\left\{ \e^{ \eta\sqrt{2\kappa} g_i\xi_i}\right\}\Big\}\\
=\bE_\eta\Big\{\prod_{i=1}^{\bar{n}}\cosh( \eta\sqrt{2\kappa} g_i)\Big\}\leq \bE_\eta\Big\{\prod_{i=1}^{\bar{n}}\e^{\kappa\eta^2 g^2_i}\Big\}=\bE_\eta\Big\{\exp\left\{\kappa\eta^2 g^Tg\right\}\Big\}=(1-2\kappa)^{-1/2}.
\end{array}
$$
Therefore,
\[
\ln \bE\Big\{\exp\big\{\kappa(\underbrace{\xi^TS\xi-1}_{\zeta})\big\}\Big\}\leq-\kappa-\half \ln(1-2\kappa)\leq {\kappa^2\over 1-2\kappa}.
\]
Now, by the standard reasoning, when optimizing with respect to $\kappa$, $0\leq \kappa<\half$ (cf. \cite[Lemma 8]{BM1998}), one conclude that
\[
\Prob\left\{\zeta\geq 2x+2\sqrt{x}\right\}\leq e^{-x},
\]
what is \rf{appeq1001}. \hfill\qed
}
\paragraph{2$^0$} The right inequality in (\ref{appeq3}) has already been justified.
To prove the left inequality in (\ref{appeq3}), we, similarly to what was done in section \ref{mainobs},
introduce the conic problem
\begin{equation}\label{appeq10}
\Opt_*=\max_{Q,t}\left\{\Tr(P^TCPQ): \;Q\succeq 0, \,\Tr(QS_k)\leq t_k\,\forall k\leq K,\, [t;1]\in\bT\right\},
\end{equation}
and  acting exactly as in the derivation of (\ref{miracle}), we arrive at
\begin{equation}\label{equal2}
\Opt=\Opt_*.
\end{equation}
Indeed, (\ref{appeq10}) is a strictly feasible and bounded conic problem, so that  its optimal value is equal to the one in its conic dual, that is,
\bse
\Opt_*&=&\min_{\lambda,[g;s],L}\left\{s: \begin{array}{l}\Tr([\sum_k\lambda_kS_k-L]Q-\sum_k[\lambda_k+g_k]t_k=\Tr(P^TCPQ)\;\;\forall (Q,t),\\
\lambda\geq0,L\succeq0,s\geq \phi_\T(-g)\\
\end{array}\right\}\\
&=&\min_{\lambda,[g;s],L}\left\{s: \begin{array}{l}\sum_k\lambda_kS_k-L=P^TCP,\,g=-\lambda,\\
\lambda\geq0,L\succeq0,s\geq \phi_\T(-g)\\
\end{array}\right\}\\
&=&\min_{\lambda}\left\{\phi_\T(\lambda): \sum_k\lambda_k S_k \succeq P^TCP, \lambda\geq0\right\}=\Opt.\\
\ese
\paragraph{3$^0$} With Lemma \ref{largedev} and (\ref{equal2}) at our disposal, we can now complete the proof of Proposition \ref{propbyproduct}
by adjusting the technique from \cite{NRT}. Specifically, problem
(\ref{appeq10}) clearly is solvable; let $Q_*,t^*$ be an optimal solution to the problem. Next, let us set $R_*=Q_*^{1/2}$, $\bar{C}=R_*P^TCPR_*$, let $\bar{C}=UDU^T$ be the eigenvalue decomposition of $\bar{C}$, and let $\bar{S}_k=U^TR_*S_kR_*U$. Observe that
\bse
\Tr(D)&=&\Tr(R_*P^TCPR_*)=\Tr(Q_*P^TCP)
=\Opt_*=\Opt,\\
\Tr(\bar{S}_k)&=&\Tr(R_*S_kR_*)=\Tr(Q_*S_k)
\leq t^*_k.
\ese
Now  let $\xi$ be Rademacher random vector.
For $k$ with $t^*_k>0$, applying Lemma \ref{largedev} to matrices $\bar{S}_k/t^*_k$, we get for $s>0$
\begin{equation}\label{appeq20}
{
\Prob\{\xi^T\bar{S}_k\xi\geq (1+2s+2\sqrt{s}) t_k^*\}\leq \e^{-s}};
\end{equation}
if $k$ is such that   $t^*_k=0$, we have $\Tr(\bar{S}_k)=0$, that is, $\bar{S}_k=0$, and (\ref{appeq20}) holds true as well.
{Observe that \rf{appeq20} implies that} that there exists a realization $\bar{\xi}$ of $\xi$ such that
$$
\bar{\xi}^T\bar{S}_k\bar{\xi}\leq {(1+2s_*+2\sqrt{s_*})}t_k^*\,\forall k
$$
{for any $s_*>\ln K$.}
Let us set $\bar{y}={1\over \sqrt{s_*}}R_*U\bar{\xi}$. Then
$$
\bar{y}^TS_k\bar{y}=s_*^{-1}\bar{\xi}^TU^TR_*S_kR_*U\bar{\xi}=s_*^{-1}\bar{\xi}^T\bar{S}_k\bar{\xi}\leq t^*_k\;\;\forall k
$$
implying that $\bar{y}\in \brX$, and
$$
\bar{y}^TP^TCP\bar{y}=s_*^{-1}\bar{\xi}^TU^TR_*CR_*U\bar{\xi}=s_*^{-1}\bar{\xi}^TD\bar{\xi}=s_*^{-1}\Tr(D)=s_*^{-1}\Opt.
$$
Thus, $\max_{y\in \brX} y^TP^TCPy\geq s_*^{-1}\Opt$, which is the first inequality in (\ref{appeq3}). \hfill\qed
%
\aic{
\subsubsection{Proof of Lemma \ref{OptSo}}\label{prpprOptSo}

\paragraph{1$^o$}  We claim that (\ref{eq:diez}) is a strictly feasible conic problem with bounded level sets of the objective (the sets where the objective is $\geq a$, for every fixed $a\in\bR$); in particular, the problem is solvable.
\par
Indeed, strict feasibility follows from the fact that the interior of the cone $\bT$ contains a positive vector, see assumptions on $\T$ in Section  \ref{sitgoal}.  Further, the projections of the feasible set onto the $[v;s]$- and $W$-spaces are bounded (the first -- since at a feasible solution it holds $0\leq s\leq1$, and the second -- due to the boundedness of the set of $v$-components of feasible solutions combined with $\sum_kS_k\succ0$). Boundedness of a level set of the objective follows from the fact that if a sequence of feasible solutions
$\{(W_i,G_i,[v^i;s^i]),i=1,2,...\}$ goes to $\infty$, then, by the above, the sequence $\{W_i,[v^i;s^i]\}$ is bounded, so that
 $\|G_i\|\to\infty$ as $i\to\infty$; since $G_i\succeq0$ due to the constraints of the problem, we have $\Tr(G_i)\to\infty$ as $i\to\infty$, which combines with boundedness of $\{W_i\}$ to imply that the objective along our sequence of feasible solutions goes to $-\infty$, which is impossible for a sequence of feasible solutions from a level set of the objective.

\paragraph{2$^o$} Our next claim is that at an optimal solution $(W,G,[v;s])$ to \rf{eq:diez} one has $s>0$.
\par Indeed, otherwise $v=0$ due to $[v;s]\in\bT$ and the origin of $\bT$, whence $W=0$ due to $W\succeq0$ and $\sum_kS_k\succ0$; besides this, $G\succeq0$, so that assuming $s=0$, we see that $\Opt_*=0$, which clearly is not the case: $\bT$ contains a vector $[\bar{v};\bar{s}]$ with, say, $\bar{s}=0.1$ and positive $\bar{v}$, implying that for some $\bar{\tau}>0$ and all $\tau\in[0,\bar{\tau}]$ tuples
$$
W_\tau=\tau I,G_\tau=[\sigma^2\bar{s}]^{-1}[BW_\tau A^TAW_\tau B^T]=[\sigma^2\bar{s}]^{-1}\tau^2BA^TAB^T,[\bar{v};\bar{s}]
$$
are feasible solutions to (\ref{eq:diez}); since $B\neq0$, for small positive $\tau$ the value of the objective of (\ref{eq:diez}) at such a solution is positive, which would be impossible when $\Opt_*=0$.
\par
{Furthermore,}
 observe that if  $(W,G,v,s)$ is an optimal solution to \rf{eq:diez} (whence, as we already know, $s>0$), when replacing $G$ with the matrix
$$
\bar{G}:=BWA^T(\sigma^2s I_m +A WA^T)^{-1}A WB^T
$$
(so that $G\succeq \bar{G}$ and $(W,\bar{G},t,s)$ is feasible for \rf{eq:diez}),
we keep the solution optimal, thus
\[
\Opt_*=\Tr\left(B [W-WA^T(\sigma^2s I_m +A WA^T)^{-1}A W]B^T\right).
\]

\paragraph{3$^o$} {To complete the proof of the lemma it suffices to show}
that the conic dual to \rf{eq:diez} is equivalent to (\ref{designproblem}); since \rf{eq:diez}, as we have already mentioned, is strictly feasible and bounded, this would imply that
$\Opt=\Opt_*.$
\par
To build the problem dual to \rf{eq:diez}, let the Lagrange
multipliers for the constraints be, respectively, $\left[\begin{array}{cc}U&V\cr V^T&Z\cr\end{array}\right]\succeq0$, $L\succeq0$, $-\lambda$, $\lambda\in\bR^K_+$, $-\tau$, $\tau\geq0$, and $[g;r]\in\bT_*$,
where
$$
\bT_*=\{[g;r]:\; r\geq \phi_{\T}(-g)\}
$$
is the cone dual to $\bT$. Taking inner products of the constraints of \rf{eq:diez} with the multipliers and summing up the results, we arrive at the aggregated constraint
$$
\begin{array}{l}
\Tr(GU) +\Tr(W[A^TV^TB+B^TVA+A^TZA+L-\sum_k\lambda_k S_k-\tau S])\\
\multicolumn{1}{r}{+\sum_k[\lambda_k+g_k]v_k+s[\sigma^2\Tr(Z)-\tau+r]+\tau\geq0}\\
\end{array}
$$
To get the dual problem, we impose on the multipliers the restriction for the resulting inequality to have the homogeneous in $W,G,v,s$ component identically equal to {\sl minus} the objective of \rf{eq:diez}, which amounts  to the relations
\bse
&&U=I_\nu,\; \tau=r+\sigma^2\Tr(Z),\;g_k=-\lambda_k\;\forall k,\\
&&[A^TV^TB+B^TVA+A^TZA+L-\sum_k\lambda_k S_k-\tau S]=-B^TB.
\ese
Under these relations, the aggregated constraint reads
$$
\Tr(BWB^T-G)\leq \tau
$$
for all feasible solutions to \rf{eq:diez}, thus $\Opt_*\leq\tau$. Therefore, the problem dual to \rf{eq:diez} is to minimize the resulting upper bound on $\Opt_*$, that is, the dual is
$$
\min_{\tau,V,Z,L,\lambda,[g;r]}\left\{\tau:\begin{array}{l}\left[\begin{array}{cc}
I_\nu&V\cr V^T&Z\cr\end{array}\right]\succeq0,\,
L\succeq0,\,\lambda\geq0,\,\tau\geq0,\,r\geq\phi_{\T}(-g)\,\\
B^TB+A^TV^TB+B^TVA+A^TZA=\sum_k\lambda_k S_k+\tau S-L\\
g=-\lambda,\tau=r+\sigma^2\Tr(Z)\\
\end{array}\right\}.
$$
Now partial minimization in $Z$ and $r$ results in $Z=V^TV$, $r=\phi_\T(-g)$, which, after eliminating $L$ and $[g;r]$, reduces the dual problem to
$$
\min_{\tau,V,\lambda}\left\{\tau:\begin{array}{l}(B+VA)^T(B+VA)\preceq\sum_k\lambda_k S_k+\tau S,\\
\lambda\geq0,\,\tau\geq \phi_{\T}(\lambda)+\sigma^2 \Tr(V^TV)\\
\end{array}\right\}.
$$
The resulting problem clearly is equivalent to (\ref{designproblem}) (substitute $V=-H^T$). Thus, (\ref{equal}) is proved.\hfill\qed

\subsubsection{Proof of Proposition \ref{prop6666}}
Under the premise of the proposition, the feasible set of (\ref{designproblementirespace}) is nonempty, and the objective clearly goes to $\infty$ along every going to $\infty$ sequence of feasible solutions $(\tau_i,H_i)$, implying that the problem is solvable. The optimal value $\Opt$ in the problem clearly is positive due to $\sigma>0$ and $B\neq0$. Now assume that (\ref{minmaxoptimality}) does not hold, so that there exists $\alpha$ and estimate $\wh{x}_*(\cdot)$ such that
\begin{equation}\label{ppcontradiction}
\alpha<\Opt\; \&\; \bE_{\xi\sim\N(0,I_m)}\{\|\wh{x}_*(Ax+\sigma\xi)-Bx\|_2^2\}\leq \alpha(1+x^TSx)\,\,\forall x\in\bR^n,
\end{equation}
and let us lead this assumption to contradiction.
\par
Consider the conic problem (cf. (\ref{eq:diez}))
\begin{equation}\label{ppeq:diezunb}
\Opt_*=\max_{W,G,s}\left\{\Tr(BWB^T)-\Tr(G):\begin{array}{l}
\left[\begin{array}{cc}G&BWA^T\cr A WB^T&\sigma^2sI_m+A WA^T\cr\end{array}\right]\succeq0,\\
W\succeq0,\,\Tr(WS)+s\leq1,\,s\geq0\\
\end{array}\right\}.
\end{equation}
This conic problem clearly is strictly feasible; the same argument as in the case of (\ref{eq:diez}) shows that the conic dual of this problem is equivalent to (\ref{designproblementirespace}) and therefore is feasible. By Conic Duality Theorem, it follows that both (\ref{ppeq:diezunb}) and (\ref{designproblementirespace}) have equal optimal values, and since $\sigma>0$, $B\neq0$, $\Opt$ is positive. Thus,
\[
\Opt_*=\Opt>0.
\]
This relation, due to $\alpha<\Opt$, implies that there is a feasible solution to (\ref{ppeq:diezunb}) with the value of the objective $>\alpha$. Since the problem is strictly feasible, feasible solutions with $s>0$ are dense in the feasible set, implying that the above feasible solution, let it be $(\widehat{W},G,\widehat{s})$, can be selected to have $\widehat{s}>0$. Further, keeping $\widehat{W}$ and $\widehat{s}$ intact and replacing $G$ with $\widehat{G}=B\widehat{W}A^T[\sigma^2\widehat{s}I_m+A \widehat{W}A^T]^{-1}A\widehat{W}B^T$, we preserve feasibility and can only increase the objective of (\ref{ppeq:diezunb}). The bottom line is that we can point out a feasible solution $(\widehat{W},\widehat{G},\widehat{s})$ to (\ref{ppeq:diezunb}) such that
\begin{equation}\label{ppgoodsolution}
\begin{array}{l}
\widehat{\alpha}:=\Tr(B^T[\widehat{W}-\widehat{W}A^T[\sigma^2\widehat{s}I_m+A \widehat{W}A^T]^{-1}A\widehat{W}]B) >\alpha,\\
\widehat{s}>0,\; \widehat{W}\succeq 0,\;\Tr(\widehat{W}S)+\widehat{s}\leq 1.
\end{array}
\end{equation}
Observe that
\begin{equation}\label{ppvarphiagain}
\widehat{\alpha}=\widehat{s}\varphi(\widehat{s}^{-1}\widehat{W})
\end{equation}
(see (\ref{diezrho})). Now let $\eta\sim \N(0,\widehat{s}^{-1}\widehat{W})$ be independent of $\xi\sim\N(0,I_m)$. We have
$$
\begin{array}{rcl}
\bE_{[\eta;\xi]}\{\|\wh{x}_*(A\eta+\sigma\xi)-B\eta\|_2^2\}&=&\bE_{\eta}\left\{\bE_{\xi}\{\|\wh{x}_*(A\eta+\sigma\xi)-B\eta\|_2^2\}\right\}\\
&\leq& \bE_{\eta}\left\{\alpha(1+\eta^TS\eta)\right\}\hbox{\ [by (\ref{ppcontradiction})]}\\
&=&\alpha(1+\widehat{s}^{-1}\Tr(\widehat{W}S)).
\end{array}
$$
By (\ref{blueeq:gaussopt}), the initial quantity in this chain is $\geq\varphi(\widehat{s}^{-1}\widehat{W})=\widehat{s}^{-1}\widehat{\alpha}$ (see (\ref{ppvarphiagain})), so that the chain yields $\widehat{s}^{-1}\widehat{\alpha}\leq \alpha(1+\widehat{s}^{-1}\Tr(\widehat{W}S))$, that is,
$$
\widehat{\alpha}\leq \alpha(\widehat{s}+\Tr(\widehat{W}S))\leq\alpha,
$$
where the last $\leq$ stems from the last inequality in (\ref{ppgoodsolution}). The resulting inequality contradicts the first inequality in (\ref{ppgoodsolution}); we have arrived at the desired contradiction.
\hfill\qed

\subsubsection{Proof of Proposition \ref{propbyproduct}}
We need the following
{\begin{lemma}\label{pplargedev} Let $S$ be a positive semidefinite $\brn\times \brn$ matrix with unit trace and $\xi$ be a Rademacher $\brn$-dimensional random vector (i.e., the entries in $\xi$ are independent and take values $\pm1$ with probabilities $1/2$). Then
\begin{equation}\label{ppappeq1}
\bE\Big\{\exp\big\{\four\xi^TS\xi\big\}\Big\}\leq 3\sqrt{2},
\end{equation}
implying that
$$
\Prob\{\xi^TS\xi>s\}\leq 3\sqrt{2}\exp\{-s/4\},\,s\geq 0.
$$
\end{lemma}}
{\bf Proof.} Let $S=\sum_{i=1}^{\brn}\lambda_ig_ig_i^T$ be the eigenvalue decomposition of $S$, so that $\lambda_i\geq0$, $\sum_i\lambda_i=1$ and $\|g_i\|_2=1$. Then
$$
\bE\Big\{\exp\big\{\four\xi^TS\xi\big\}\Big\}=\bE\{\exp\{{1\over 4}\sum_i\lambda_i(g_i^T\xi)^2\}\}
$$
is a convex function of $\lambda$ and therefore it attains its maximum over nonnegative vectors $\lambda$ with unit sum of entries at a basic orth. Thus, it suffices to verify (\ref{ppappeq1}) when $S=gg^T$ with unit vector $g$. {By the Hoeffding inequality one has
\[
\Prob\{|g^T\xi|>s\}\leq 2\exp\{-s^2/2\}.
\]}
It follows that $\Prob\{(g^T\xi)^2>r\}\leq 2\exp\{-r/2\}$, and thus $p(r):=\Prob\{{1\over 4}(g^T\xi)^2\geq r\}\leq 2\exp\{-2r\}$. Consequently,
$$
\begin{array}{l}
\bE\Big\{\exp\{\four(g^T\xi)^2\}\Big\}=\int_0^\infty \exp\{r\}[-dp(r)] =\int_0^\infty \exp\{r\}p(r)dr+1\\
\leq
\int_0^\infty \exp\{r\}\min[1,2\exp\{-2r\}]dr + 1=1+\int_0^{{1\over 2}\ln(2)}\exp\{r\}dr +2\int_{{1\over 2}\ln(2)}^\infty \exp\{-r\}dr=3\sqrt{2}.
\end{array}\eqno{\hbox{\qed}}
$$

\paragraph{2$^0$} The right inequality in (\ref{appeq3}) has already been justified.
To prove the left inequality in (\ref{appeq3}), we, similarly to what was done in Section \ref{mainobs},
introduce the conic problem
\begin{equation}\label{ppappeq10}
\Opt_*=\max_{Q,t}\left\{\Tr(P^TCPQ): Q\succeq 0, \Tr(QS_k)\leq t_k\,\forall k\leq K, [t;1]\in\bT\right\},
\end{equation}
and  acting exactly as in the derivation of (\ref{miracle}), we arrive at
\begin{equation}\label{ppequal2}
\Opt=\Opt_*.
\end{equation}
Indeed, (\ref{ppappeq10}) is a strictly feasible and bounded conic problem, so that  its optimal value is equal to the one in its conic dual, that is,
\bse
\Opt_*&=&\min_{\lambda,[g;s],L}\left\{s: \begin{array}{l}\Tr([\sum_k\lambda_kS_k-L]Q)-\sum_k[\lambda_k+g_k]t_k=\Tr(P^TCPQ)\;\;\forall (Q,t),\\
\lambda\geq0,L\succeq0,s\geq \phi_\T(-g)\\
\end{array}\right\}\\
&=&\min_{\lambda,[g;s],L}\left\{s: \begin{array}{l}\sum_k\lambda_kS_k-L=P^TCP,\,g=-\lambda,\\
\lambda\geq0,L\succeq0,s\geq \phi_\T(-g)\\
\end{array}\right\}\\
&=&\min_{\lambda}\left\{\phi_\T(\lambda): \sum_k\lambda_k S_k \succeq P^TCP, \lambda\geq0\right\}=\Opt.\\
\ese

\paragraph{3$^0$} With Lemma \ref{pplargedev} and (\ref{ppequal2}) at our disposal, we can now complete the proof of Proposition \ref{propbyproduct}
by adjusting the technique from \cite{NRT}. Specifically, problem
(\ref{ppappeq10}) clearly is solvable; let $Q_*,t^*$ be an optimal solution to the problem. Next, let us set $R_*=Q_*^{1/2}$, $\bar{C}=R_*P^TCPR_*$, let $\bar{C}=UDU^T$ be the eigenvalue decomposition of $\bar{C}$, and let $\bar{S}_k=U^TR_*S_kR_*U$. Observe that
\bse
\Tr(D)&=&\Tr(R_*P^TCPR_*)=\Tr(Q_*P^TCP)
=\Opt_*=\Opt,\\
\Tr(\bar{S}_k)&=&\Tr(R_*S_kR_*)=\Tr(Q_*S_k)
\leq t^*_k.
\ese
Now  let $\xi$ be Rademacher random vector.
For $k$ with $t^*_k>0$, applying Lemma \ref{pplargedev} to matrices $\bar{S}_k/t^*_k$, we get for $s>0$
\begin{equation}\label{ppappeq20}
\Prob\{\xi^T\bar{S}_k\xi> s t_k^*\}\leq 3\sqrt{2}\exp\{-s/4\};
\end{equation}
if $k$ is such that   $t^*_k=0$, we have $\Tr(\bar{S}_k)=0$, that is, $\bar{S}_k=0$, and (\ref{ppappeq20}) holds true as well.
Now let
$$
s_*=4\ln(5K),
$$
so that $3\sqrt{2}\exp\{-s_*/4\}<1/K$. The latter relation combines with (\ref{ppappeq20}) to imply that there exists a realization $\bar{\xi}$ of $\xi$ such that
$$
\bar{\xi}^T\bar{S}_k\bar{\xi}\leq s_*t_k^*\,\forall k.
$$
Let us set $\bar{y}={1\over \sqrt{s_*}}R_*U\bar{\xi}$. Then
$$
\bar{y}^TS_k\bar{y}=s_*^{-1}\bar{\xi}^TU^TR_*S_kR_*U\bar{\xi}=s_*^{-1}\bar{\xi}^T\bar{S}_k\bar{\xi}\leq t^*_k\;\;\forall k
$$
implying that $\bar{y}\in \brX$, and
$$
\bar{y}^TP^TCP\bar{y}=s_*^{-1}\bar{\xi}^TU^TR_*CR_*U\bar{\xi}=s_*^{-1}\bar{\xi}^TD\bar{\xi}=s_*^{-1}\Tr(D)=s_*^{-1}\Opt.
$$
Thus, $\max_{y\in \brX} y^TP^TCPy\geq s_*^{-1}\Opt$, which is the first inequality in (\ref{appeq3}). \hfill\qed
}{}
\section{Calculus of Ellitopes}\label{ppapp:ovalcalc}
\begin{itemize}
\item Intersection ${\X}=\bigcap\limits_{i=1}^I {\X}_i$ of {ellitopes} ${\X}_i=\{x\in\bR^n: \exists (y^i\in\bR^{n_i}, t^i\in\T_i): x=P_iy^i\ \&\ [y^i]^TS_{ik}y^i\leq t^i_k,1\leq k\leq K_i\}$, is an {ellitope}. Indeed, this is evident when $\X=\{0\}$. Assuming $\X\neq\{0\}$, we have
$$
\begin{array}{rcl}
{\X}&=&\{x\in\bR^n:\exists (y=[y^1;...;y^I]\in \Y, t=(t^1,...,t^I)\in\T=\T_1\times...\times \T_I):\\
 &&\multicolumn{1}{r}{x=Py:=P_1y^1\ \&\ \underbrace{[y^i]^TS_{ik}y^i}_{y^TS^+_{ik}y}\leq t^i_k,1\leq k\leq K_i,1\leq i\leq I\},}\\
\Y&=&\{[y^1;...;y^I]\in\bR^{n_1+...+n_I}: P_iy^i=P_1y^1,\,2\leq i\leq I\}\\
\end{array}
$$
(note that $\Y$ can be identified with $\bR^{\brn}$ with a properly selected $\brn>0$);
\item Direct product ${\X}=\prod\limits_{i=1}^I {\X}_i$  of {ellitopes} ${\X}_i=\{x^i\in\bR^{n_i}: \exists (y^i\in\bR^{\brn_i},t^i\in\T_i): x^i=P_iy^i,\,1\leq i\leq I\ \&\ [y^i]^TS_{ik}y^i\leq t^i_k,1\leq k\leq K_i\}$ is an {ellitope}:
$$
\begin{array}{l}
{\X}=\{[x^1;...;x^I]\in\bR^{n_1}\times...\times \bR^{n_I}:\exists \left(\begin{array}{c}y=[y^1;...;y^I]\in\bR^{\brn_1+...\brn_I}\\
t=(t^1,...,t^I)\in\T=\T_1\times...\times \T_I\\
\end{array}\right))\\
 \multicolumn{1}{r}{
 x=Py:=[P_1y^1;...;P_Iy^I],\,
 \underbrace{[y^i]^TS_{ik}y^i}_{y^TS^+_{ik}y}\leq t^i_k,1\leq k\leq K_i,1\leq i\leq I\}}\\
 \end{array}
$$
\item The linear image ${\Z}=\{Rx:x\in {\X}\}$, $R\in\bR^{p\times n}$,  of an {ellitope} ${\X}=\{x\in\bR^n:\exists (y\in\bR^{\brn},t\in\T): x=P_y\ \&\ y^TS_ky\leq t_k,\,1\leq k\leq K\}$ is an {ellitope}:
$$
{\Z}=\{z\in \bR^p: \exists (y\in\bR^{\brn},t\in\T): z=[RP]y\ \&\ y^TS_ky\leq t_k,\,1\leq k\leq K\}.
$$
\item The inverse linear image ${\Z}=\{z\in \bR^q: Rz\in {\X}\}$, $R\in\bR^{n\times q}$, of an {ellitope} ${\X}=\{x\in\bR^n: \exists (y\in\bR^{\brn},t\in\T): x=Py\ \&\ y^TS_ky\leq t_k,1\leq k\leq K\}$ under linear mapping $z\mapsto Rz:\bR^q\to\bR^n$ is an {ellitope}, {\sl provided that the mapping is an embedding: $\Ker R=\{0\}$}:
$$
\begin{array}{rcl}
{\Z}&=&\{z\in\bR^q:\exists (y\in \Y,t\in\T): z=\bar{P}y\ \&\ y^TS_ky\leq t_k,\,1\leq k\leq K\},\\
\Y&=&\{y\in\bR^{\brn}:Py\in\hbox{\rm Im} R\},\\
\bar{P}&:& \bar{P}y = \Pi R,\,\hbox{\ where $\Pi:\hbox{\rm Im} R\to \bR^q$ is the inverse of $z\mapsto Rz: \bR^q\to \hbox{\rm Im} R$}\\
\end{array}
$$
($\Y$ can be identified with some $\bR^k$, and $\Pi$ is well defined since $R$ is an embedding).
\item The arithmetic sum ${\X}=\{x=\sum_{i=1}^I x^i:x^i\in {\X}_i,\,1\leq i\leq I\}$, of {ellitopes} ${\X}_i$ is an {ellitope}, with representation readily given by those of ${\X}_1,...,{\X}_I$. \par
Indeed, ${\X}$ is the image of ${\X}_1\times,,,\times {\X}_I$ under the linear mapping $[x^1;...;x^I]\mapsto x^1+....+x^I$, and taking direct products and images under linear mappings preserve {ellitopes}.

\end{itemize}
 Note that the outlined ``calculus rules'' are fully algorithmic: representation (\ref{setX}) of the result of an operation is readily given by the representations (\ref{setX}) of the operands.

\section{Numerical lower bounds of the minimax risk}
\label{app:lbc}
To implement efficiently the bounding scheme sketched in section \ref{sec:num1} we need to provide a convex (and numerically tractable) set $\Q_\delta$ of covariance matrices $Q$ such that for any $Q\in \Q_\delta$,
\[
\Prob_{\eta\sim\N(0,Q)}\{\eta\in \X\}=\Prob_{\eta\sim\N(0,Q)}\left\{\exists t\in\T:\;\eta^TS_k\eta\leq t_k,\,1\leq k\leq K,\right\}\geq 1-\delta.
\]
Such sets can be constructed straightforwardly in the case where $\X$ is an ellipsoid or a parallelotope (e.g., a box).
\par
To build a lower bound for an optimal risk on the ellipsoid
\[
\X_1=\{x\in \bR^n:\;x^TS_1x\leq 1\},\;\;
\]
where $S_1\succ 0$ is a given matrix, recall that for any $\beta>2\max_i w_{i}$, where $w_{i}$ are the eigenvalues of the matrix $W=[S_1^{1/2}QS_1^{1/2}]$ we have (see, e.g., Lemma \ref{lem1})
\be
\Prob_{\eta\sim \N(0,Q)}\{\eta^TS_1\eta>1\}&\leq&\exp\left\{-\half\sum\limits_{i=1}^n\ln(1-2\beta^{-1}w_{i})- \beta^{-1}\right\}\nn
&=& \exp\left\{-\half\ln\,\Det\left[I-2\beta^{-1}{W}\right]-\beta^{-1}\right\}
\ee{eq:ex-deltarho}
(this relation is given by the first equality in \rf{then1} with $S_1$ in the role of $S$ and $\beta=1/\gamma$).
Let now $\delta>0$; we conclude that for all $Q\in \Q_{1,\delta}$ where
\[
\begin{array}{rl}
\Q_{1,\delta}=\Big\{Q\succeq0:\,&\exists \beta>0:\,\beta I\succ 2W,\,W=[S_1^{1/2}QS_1^{1/2}]\\
&-{\beta\over 2}\ln\,\Det\left[I-2\beta^{-1}{W}\right]+\beta\ln(1/\delta)\leq 1,\Big\},
\end{array}
\]
one has $\Prob_{\eta\sim \N(0,Q)}\{\eta\notin \X\}\leq \delta$.
Though efficiently tractable, the set ${\Q}_{1,\delta}$ is still difficult to deal with numerically -- solving the problem
\be
\min_{Q\in \Q_{1,\delta}}\varphi(Q)
\ee{eq:optd}
(e.g., using CVX) takes hours already for small problem sizes. Therefore, in the experiments presented in section \ref{sec:num1} we used two simple substitutes
\begin{itemize}
\item[[1.]] an appropriate ``contraction'' $\Q_{\rho,\delta}$ of $\Q:=\left\{Q\succeq 0:\,\Tr(QS_1)\leq 1\right\}$:
\[
\Q_{\rho,\delta}=\{Q\succeq 0:\,\Tr(QS_1)\leq \rho\},
\]
where $\rho$ was chosen according to Lemma \ref{lem1} to ensure that $\Prob_{\eta\sim \N(0,Q)}\{\eta\notin \X_1\}\leq \delta$ for all $Q\in \Q_{\rho,\delta}$. This construction underlies the lower bound represented by red curves on Figures \ref{fig:e1-1} and \ref{fig:e1-end};
\item[[2.]] a ``quadratic approximation'' $\overline{\Q}_{1,\delta}\subseteq \Q_{1,\delta}$ (see, e.g., \cite[Lemma 1]{laurent2000adaptive}):
\[
\overline{\Q}_{1,\delta}=\left\{Q\succeq 0:\,\Tr(QS_1)+2\|QS_1\|_2\sqrt{\ln(1/\delta)}+2\|QS_1\|\ln(1/\delta)\leq 1\right\}.
\]
This approximation is used to compute the lower bound represented by magenta curves on Figures \ref{fig:e1-1} and \ref{fig:e1-end}.
\end{itemize}
The strategy we use to compute the lower estimates for $\Riskopt[\X_1]$ amounts to solve the problem \rf{eq:optd} with ${\Q}_{\rho,\delta}$ or $ \overline{\Q}_{1,\delta}$ in the role of ${\Q}_\delta$ for several values of $\delta$. Then the bound  \rf{eq:ex-deltarho} with the computed optimal solution $Q$ to \rf{eq:optd} is used to obtain a refined estimate $\delta'$ for the probability that $\eta\not \in \X$ which is then substituted into the lower risk bound \rf{eq:lrb}; due to its origin, $\delta'\leq\delta$.
\par
Let us now consider the situation where the set of signals is a parallelotope, namely,
\[
\X_2=\{x\in \bR^n:\;x^TS_kx\leq 1,\,S_k=a_ka_k^T,\, a_k\in\bR^n, 1\leq k\leq K\}.
\]
When $\eta\sim \N(0,Q)$, random variables $\zeta_k=a_k^T\eta$ follow Gaussian distribution, $\zeta_k\sim \N(0, a_k^TQa_k)$,  and $\Prob\{\zeta_k^2\geq 1\}\leq \alpha$, $0< \alpha\leq 1$ if $a_k^TQa_kq^{2}_{1-\alpha/2}\leq 1$ where $q_\beta$ is the $\beta$-quantile of the standard normal distribution.
We conclude that if $Q\in \Q_{2,\delta}$,
\[
Q_{2,\delta}:=\big\{Q\succeq 0,\,a_k^TQa_kq^{2}_{1-\delta/(2K)}\leq 1, \,1\leq k\leq K\big\},
\]
one has $\Prob_{\eta\sim \N(0,Q)}\{\eta\notin \X_2\}\geq 1-\delta$. Finally, to lower bound $\Riskopt[\X_2]$, given $0<\delta<1$, we compute an optimal solution to \rf{eq:optd} with $\Q_{1,\delta}$ replaced with $\Q_{2,\delta}$, and then apply the bound \rf{eq:0002}.
\section{Conic duality}\label{conicd}
A conic problem is an optimization problem of the form
$$
\Opt(P)=\max_x\left\{c^Tx: A_ix-b_i\in\bK_i,i=1,...,m, Px=p\right\}\eqno{(P)}
$$
where $\bK_i$ are regular (i.e., closed, convex, pointed and with a nonempty interior) cones in Euclidean spaces $E_i$. Conic dual of $(P)$ ``is responsible'' for upper-bounding the optimal value in $(P)$ and is built as follows: selecting somehow  {\sl Lagrange multipliers} $\lambda_i$ for the conic constraints $A_ix-b_i\in\bK_i$ in the cones dual to $\bK_i$:
$$
\lambda_i\in \bK_i^*:=\{\lambda:\langle\lambda,y\rangle\geq0\,\forall y\in\bK_i\},
$$
and a Lagrange multiplier $\mu\in\bR^{\dim p}$ for the equality constraints,
every feasible solution $x$ to $(P)$ satisfies the linear inequalities $\langle\lambda_i,A_ix\rangle\geq \langle \lambda_i,b_i\rangle$, $i\leq m$, same as the inequality $\mu^TPx\geq \mu^Tp$, and thus satisfies the aggregated inequality
$$
\sum_i\langle\lambda_i,A_ix\rangle +\mu^TPx\geq \sum_i\langle \lambda_i,b_i\rangle +\mu^Tp.
$$
If the left hand side of this inequality is, {\sl identically in $x$}, equal to $-c^Tx$ (or, which is the same, $-c=\sum_iA_i^*\lambda_i+P^T\mu$, where $A_i^*$ is the conjugate of $A_i$), the inequality produces an upper bound $-\langle \lambda_i,b_i\rangle-p^T\mu$ on $\Opt(P)$. The dual problem
$$
\Opt(D)=\min_{\lambda_1,...,\lambda_m,\mu}\left\{-\sum_i\langle \lambda_i,b_i\rangle-p^T\mu:\lambda_i\in\bK_i^*,i\leq m, \sum_iA_i^*\lambda_i+P^T\mu=-c\right\}\eqno{(D)}
$$
is the problem of minimizing this upper bound. Note that $(D)$ is a conic problem along with $(P)$ -- it is a problem of optimizing a linear objective under a bunch of linear equality constraints and conic inclusions of the form ``affine function of the decision vector should belong to a given regular cone.'' Conic Duality Theorem (see, e.g., \cite{LMCO}) states that when one of the problems $(P)$, $(D)$ is bounded\footnote{for a maximization (minimization)  problem, boundedness means that the objective is bounded from above (resp., from below) on the feasible set.} and strictly feasible, then the other problem in the pair is solvable, and $\Opt(P)=\Opt(D)$. In this context, strict feasibility exactly means that there exists a feasible solution for which all conic inclusions are satisfied strictly, that is, the left hand side of the inclusion belongs to the {\sl interior} of the right hand side cone.
%
%
%


\end{document}